\documentclass[12pt]{amsart}
\usepackage[]{hyperref} 
\usepackage{amssymb}
\usepackage[curve]{xypic}
\newcommand{\Omit}[1]{\begin{tiny}#1\end{tiny}}
\renewcommand{\Omit}[1]{}

\newbox\mybox
\def\overtag#1#2#3{\setbox\mybox\hbox{$#1$}\hbox to
  0pt{\vbox to 0pt{\vglue-#3\vglue-\ht\mybox\hbox to \wd\mybox
      {\hss$\scriptstyle#2$\hss}\vss}\hss}\box\mybox}
\def\undertag#1#2#3{\setbox\mybox\hbox{$#1$}\hbox to 0pt{\vbox to
    0pt{\vglue#3\vglue\ht\mybox\hbox to \wd\mybox
      {\hss$\scriptstyle#2$\hss}\vss}\hss}\box\mybox}
\def\lefttag#1#2#3{\hbox to 0pt{\vbox to 0pt{\vss\hbox to
      0pt{\hss$\scriptstyle#2$\hskip#3}\vss}}#1}
\def\righttag#1#2#3{\hbox to 0pt{\vbox to 0pt{\vss\hbox to
      0pt{\hskip#3$\scriptstyle#2$\hss}\vss}}#1}

\def\Dot{\lower.2pc\hbox to 2.5pt{\hss$\bullet$\hss}}
\def\Circ{\lower.2pc\hbox to 2.5pt{\hss$\circ$\hss}}
\def\Vdots{\raise5pt\hbox{$\vdots$}}
\def\splicediag#1#2{\xymatrix@R=#1pt@C=#2pt@M=0pt@W=0pt@H=0pt}

\renewcommand\frame[2][3pt]{\hbox{$\vcenter{\hbox{\vrule\vbox
{\hrule\kern#1\hbox{\kern#1$#2$\kern#1}\kern#1\hrule}\vrule}}$}}

\newcommand\lineto{\ar@{-}}
\newcommand\dashto{\ar@{--}}
\newcommand\dotto{\ar@{.}}

\newcommand{\C}{\mathbb C}

\newcommand{\Q}{\mathbb Q}

\newtheorem*{RatConj}{Rational Conjecture}

\newtheorem*{Main}{Main Conjecture of \cite{mu}}

\newtheorem*{theorem*}{Theorem}
\newtheorem{theorem}{Theorem}[section]
\newtheorem{proposition}[theorem]{Proposition}
\newtheorem*{proposition*}{Proposition}
\newtheorem{lemma}[theorem]{Lemma}
\newtheorem{corollary}[theorem]{Corollary}
\newtheorem*{corollary*}{Corollary}
\theoremstyle{definition}
\newtheorem{example}[theorem]{Example}
\newtheorem*{example*}{Example}
\newtheorem{definition}[theorem]{Definition}
\newtheorem{remark}[theorem]{Remark}
\newtheorem*{remark*}{Remark}

\newtheorem*{conjecture*}{Conjecture}
\begin{document}

\title{The number of equisingular moduli of a rational surface singularity}
\author{Jonathan Wahl}
\dedicatory{To Henry Laufer on his 70th birthday}
\address{Department of Mathematics\\The University of North
  Carolina\\Chapel Hill, NC 27599-3250} \email{jmwahl@email.unc.edu}
  \begin{abstract} 
  We consider a conjectured topological inequality for the number of equisingular moduli of a rational surface singularity, and prove it in some natural special cases. When the resolution dual graph is ``sufficiently negative'' (in a precise sense), we verify the inequality via an easy cohomological vanishing theorem, which implies that this number is computed simply from the graph (Theorem \ref{rt}).   To consider an important and less restrictive meaning of ``sufficiently negative'' requires a much more difficult ``hard vanishing theorem'' (Theorem \ref{t1}), which is false in characteristic $p$.  Theorem \ref{t3} verifies the conjectured inequality in this more general situation.  As a corollary, we classify in characteristic $p$ all taut singularities with reduced fundamental cycle (Theorem \ref{tautp}).
  \end{abstract}
\maketitle
\section{Introduction}
The following is a special case of a conjectured inequality in \cite{mu} for complex normal surface singularities:
\begin{RatConj} Let $(V,0)$ be (the germ of) a complex rational surface singularity, $(X,E)\rightarrow (V,0)$ the minimal good resolution (or MGR).  Define $S_X=(\Omega^1_X(\text{log}\ E))^*$, the bundle on $X$ of derivations logarithmic along $E$.  Then $$h^1(X,S_X)\leq h^1(X,-(K_X+E)),$$ with equality if and only if $(V,0)$ is quasi-homogeneous.
\end{RatConj}

The right-hand term is $h^1(X,\wedge^2S_X)$, the \emph{second plurigenus} of the singularity, which can be computed from the resolution graph $\Gamma$, hence is a topological invariant.  On the other hand, $h^1(S_X)$ is the dimension of the smooth space of \emph{equisingular deformations} of the  singularity (\cite{wan}, (5.16)), i.e., those deformations with the same graph; it is difficult to compute, and can vary in equisingular families.  In \cite{mu} the Conjecture was verified when the graph $\Gamma$ is star-shaped, or when $(V,0)$ admits a smoothing whose total space is $(\C^3/G,0)$.

 In this paper we prove the Conjecture for graphs $\Gamma$ which are ``sufficiently negative at the nodes'' by computing $h^1(S_X)$ in all those cases.  Let $r$ be the number of ends of the graph $\Gamma$, and for an exceptional curve $E_i$ let $t_i$ be its valence and $d_i=-E_i\cdot E_i$.  The easiest version is the following:
 
\begin{theorem*}[\ref{rt}] Let $(X,E)\rightarrow (V,0)$ be the MGR of a rational surface singularity.  Suppose that for all $i$, one has $$d_i\geq 2t_i-2.$$
\begin{enumerate}
\item If $\Gamma$ is star-shaped, then $(V,0)$ is weighted homogeneous, and 

$h^1(S_X)=h^1(-(K_X+E))=r-3.$
\item If $\Gamma$ is not star-shaped, then $h^1(S_X)=h^1(-(K_X+E))-1=r-4.$
\end{enumerate}
\end{theorem*}
  For instance, if $\Gamma$ is any trivalent tree (not star-shaped) with $d\geq 4$ at each node, then the ``number of equisingular moduli" is exactly the number of ends minus $4$.  Note however that the reduced curve $E$ is itself rigid.
  
The base space of a semi-universal deformation of $(V,0)$ contains a unique irreducible \emph{Artin} component, parametrizing deformations which resolve simultaneously after base change (see e.g. \cite{ste}, p.\ $115$).  It is smooth of  dimension $h^1(\Theta_X)=h^1(S_X)+\Sigma(d_i-1)$.  Combining with the results of \cite{dej} yields

\begin{corollary*}[\ref{rtc}] For a rational singularity $(V,0)$ with $d_i\geq 2t_i-2$, all $i$, one has
\begin{enumerate}
  \item if $\Gamma$ is star-shaped, then dim $T^1_V=\sum_i(2d_i-3)+r-4.$
     \item if $\Gamma$ is not star-shaped, then dim $T^1_V=\sum_i(2d_i-3)+r-5.$
\end{enumerate}
\end{corollary*}

 The origin of the Rational Conjecture is the
\begin{Main}  Let $(X,E)\rightarrow (V,0)$ be the MGR of a complex normal surface singularity.  Define $S_X=(\Omega^1_X(\text{log}\ E))^*$.  If $(V,0)$ is not Gorenstein, then $$h^1(\mathcal O_X)-h^1(S_X)+h^1(-(K_X+E))\geq 0,$$ with equality if and only if $(V,0)$ is quasihomogeneous.
\end{Main}

For $(V,0)$ Gorenstein, it followed from \cite{w} that the cohomological expression above is always greater than or equal to $1$, with equality exactly in the quasihomogeneous case.  That result was a generalization of an inequality for isolated hypersurface singularities: the Milnor number $\mu$ is greater than or equal to the Tjurina number $\tau$, with equality exactly in the quasihomogeneous case.  The connection is that for a two-dimensional hypersurface, on the minimal good resolution one has $$1+\mu-\tau=h^1(\mathcal O_X)-h^1(S_X)+h^1(-(K_X+E)).$$ The reader may consult \cite{mu} to see the derivation of the expression in the Main Conjecture, the verification for certain quotients of hypersurface singularities, and the proof that one does have equality in the quasihomogeneous case.

 In this more general setting, $h^1(S_X)$ counts the first-order deformations of $X$ to which each exceptional curve $E_i$ lifts.  It is the tangent space to a smooth family of equisingular deformations of the \emph{resolution}; but only sometimes (e.g., when $h^1(\mathcal O_X)=h^1(\mathcal O_E)$) do these deformations ``blow down'' to actual deformations of the singularity (\cite{wan}, (2.7)).

Our initial approach to the Rational Conjecture applies as well to the Main Conjecture.  By Propositions \ref{p6} and \ref{p2.10}, one can compute $h^1$ of the restrictions of the $3$ bundles to $E$ solely from the graph $\Gamma$. The following completely general result is instructive.

\begin{proposition*}[\ref{p11}] Let $(X,E)\rightarrow (V,0)$ be the MGR of a normal surface singularity, not a simple elliptic or cusp singularity.   Then $$h^1(\mathcal O_E)-h^1(S_X\otimes \mathcal O_E)+h^1(-(K_X+E)\otimes \mathcal O_E)=1-\delta,$$
where $\delta$ is $1$ if the resolution dual graph $\Gamma$ is star-shaped, $0$ otherwise.  
\end{proposition*}

From this result follows the key observation: the inequality in the Main Conjecture holds in those cases for which $$h^1(S_X)=h^1(S_X\otimes \mathcal O_E),$$ e.g. if $h^1(S_X(-E))=0.$  
  If the $d_i$ are big enough relative to $g_i=\text{genus} \ E_i$ and $t_i$, then such a cohomological vanishing result can be proved via so-called ``easy vanishing theorems'' of \cite{vt}, as recalled in Section $3$.   Deducing quasihomogeneity is often possible via a result in \cite{w} (Proposition \ref{w} below).

\begin{theorem*}[\ref{vtt}] Let $(X,E)\rightarrow (V,0)$ be the MGR of a normal surface singularity.  Suppose that for all $i$ one has $$d_i\geq 4g_i-4+3t_i,$$  with strict inequality for at least one $i$.  Then $h^1$ of each of the three sheaves $\mathcal O_X, S_X$, and $ -(K_X+E))$ is equal to $h^1$ of its restriction to $E$, and the Main Conjecture holds.
\end{theorem*}
This should be compared with H. Grauert's old result: if $E$ is a single smooth curve and $d>4g-4$, then $(V,0)$ is the cone over a curve, determined by $E$ and its conormal bundle.  Here, $h^1(S_X)=h^1(S_X\otimes \mathcal O_E)=(3g-3)+g.$

Returning to the Rational Conjecture, it is desirable (and necessary for applications) to prove it in a slightly more general situation than Theorem \ref{rt}, with a weaker inequality for the $d_i$.   However, this requires a much more delicate argument and  a ``hard vanishing theorem'' (as in \cite{w}), which is false in prime characteristic.  The main work of the paper is to prove vanishing under certain conditions.

Let $E'$ denote the sum of the curves which are not end-curves; and for each curve, let $t'_i$ denote the number of intersections with curves in $E'$.  The condition $$(**)\ d_i\geq t_i+t'_i-2, \text{ all i} $$ is exactly what simplifies the computation of  $h^1(-(K_X+E))$; by Proposition \ref{t'}, it equals $r-3$ (except for cyclic quotients).  As $h^1(S_X\otimes \mathcal O_E)$ is $r-4$ in the non-star-shaped case (Proposition \ref{p6}) and $h^1(S_X)\geq h^1(S_X\otimes \mathcal O_E)$, we deduce:
 \begin{proposition*}[\ref{p4}]  Let $(V,0)$ be a rational singularity whose graph satisfies $(**)$.  Then the Rational Conjecture for $(V,0)$ is equivalent to
$$H^1(S_X(-E'))=0.$$  
\end{proposition*}

 A long and technical argument over several sections of the paper eventually yields the following, which is somewhat weaker than desired:

\begin{theorem*}[\ref{t1}] Let $(V,0)$ be a rational singularity whose graph satisfies $(**)$.  Then $H^1(S_X(-E')\otimes \mathcal O_E)=0$.
\end{theorem*}
This result does yield the Rational Conjecture in an important case.
\begin{corollary*}[\ref{c1}]  Suppose a rational singularity with reduced fundamental cycle satisfies $(**)$.  Then $H^1(S_X(-E'))=0$ and $h^1(S_X)=r-4+\delta.$
\end{corollary*}
 The set-up used  above can be stretched to prove a more precise result.
\begin{theorem*}[\ref{t3}]  Suppose a rational singularity, with reduced fundamental cycle, has all curves satisfying $d_i\geq t_i+t'_i-2$ for all $i$, except that one also allows that either
\begin {enumerate}
\item exactly one curve satisfies $d=t+t'-3$, or 
\item exactly two curves, separated by a (possibly empty) string of rational curves, satisfy $d=t+t'-3$.
\end{enumerate}
Then $h^1(S_X(-E'))=0$, $h^1(S_X)=h^1(S_X\otimes  \mathcal O_E)$, and the Rational Conjecture holds.
\end{theorem*} 
(If some $d \leq t+t'-4$, then $h^1(S_X(-E')\otimes\mathcal O_E)\neq 0$ (Lemma 5.1), and $h^1(S_X)> h^1(S_X\otimes \mathcal O_E)$.)

     \begin{remark*}  
      Some of the issues in this paper originate with the work of Henry Laufer \cite{lauf}, \cite{laufert}.  In \cite{laufert}, he  made a complete list of graphs $\Gamma$ for which there corresponds a unique analytic type; these are the \emph{taut} singularities, characterized by the vanishing of $h^1(S_X)$ for \emph{every} singularity with graph $\Gamma$.  He also listed those $\Gamma$ for which the singularity is uniquely determined by the analytic type of the reduced curve $E$ (i.e.,  in the rational case, by cross-ratios at the nodes).  Theorem \ref{t3} allows one to recover easily these classifications for rational singularities with reduced fundamental cycle.   (For instance, Laufer's condition on line $5$ of p.$162$ of \cite{laufert} is equivalent to $d\geq t+t'-3$.) Of course, these form a very small portion of Laufer's lists!
      

\end{remark*}



More importantly,  the methods of this paper allow one to extend this partial classification to characteristic $p$.  In \cite{ma}, M. Artin listed all the rational double points in characteristic $p$; he showed that for a graph of type $D$ or $E$, there were certain primes for which the graph was not taut. Lee-Nakayama \cite{LN} showed that all cyclic quotients are taut.  Work of F. Sch\"{u}ller \cite{sch} extends some work of Laufer to characteristic $p$, so that a graph $\Gamma$ is taut if and only if $h^1(S_X)=0$ for every singularity with that graph.  
\begin{theorem*}[\ref{tautp}]  In characteristic $p$, there is a complete list of the taut singularities with reduced fundamental cycle.
\end{theorem*}

The paper is organized as follows: In Section $2$ we compute explicitly $h^0$ and $h^1$ of the restrictions of the three relevant sheaves to $E$, finding (Proposition \ref{p11}) an equality close to the Main Conjecture.  The divisor $E'$ (which is $E$ minus the rational end-curves) becomes important.  In Section $3$, ``easy vanishing theorems'' (in the sense of \cite{vt}) give explicit conditions, in both the general and rational cases, for $h^1$'s of the relevant bundles on $X$ to be computable from restriction to $E$.  From  Section $4$ on, only rational singularities are considered, and one attempts to get weaker restrictions on the $d_i$ to imply the vanishing of $H^1(S_X(-E'))$.  This is the technical heart of the paper: one uses an inductive procedure on subgraphs, starting at the end-curves and growing towards the interior.  Section $7$ generalizes the preceding and gives the strongest results of the paper.  Finally, in Section $8$ previous proofs are examined and modified to get analogous results valid in characteristic $p$.  The taut singularities with reduced fundamental cycle are listed in Section $9$.

\section{Restriction to $E$ and $E'$}

Consider the minimal good resolution $(X,E)\rightarrow (V,0)$ of a normal surface singularity, with weighted dual graph $\Gamma$.  For each exceptional curve $E_i$, let $g_i$ be the genus, $d_i=-E_i \cdot E_i$ the degree of the conormal bundle,
and $t_i$ the number of intersections with other curves (or, valency of the vertex in $\Gamma$).   A curve (or vertex) is called a \emph{node} if $t_i\geq 3$, and a \emph{star} if $2g_i+t_i>2.$  

   Let $S=S_X$ be the rank $2$ bundle of derivations logarithmic along $E$, defined by the short exact sequence $$0\rightarrow  S\rightarrow \Theta_X\rightarrow \oplus N_{E_{i}}\rightarrow 0.$$   In local coordinates $x,y$, if $E$ is given by $y=0$, then $S$ is generated by $\partial/\partial x$ and $y\partial/\partial y$; if $E$ is given by $xy=0$, $S$ is generated by $x\partial/\partial x$ and $y\partial/\partial y$.  For each $i$ there is an exact sequence $$0\rightarrow \mathcal O_{E_i}\rightarrow S\otimes \mathcal O_{E_i}\rightarrow \Theta_{E_i}(-t_i)\rightarrow 0.$$ (We abuse notation slightly, as $-t_i$ represents the negative of an effective divisor of degree $t_i$.)
This sequence splits unless $E$ consists of a single smooth curve.  The global section of $S\otimes \mathcal O_{E_i}$ from the left hand injection sends $1$ to $y\partial/\partial y$, where $y=0$ is \emph{any} local equation for $E_i$ on $X$; for, $y'\partial/\partial y'=y\partial/\partial y$ modulo $y=0$.   We record a useful
\begin{lemma}  If $E_i\cap E_j=P_{ij}$, then $H^0(S\otimes \mathbb C_{P_{ij}})$ has a natural ordered basis $\{x\partial/ \partial x, y\partial /\partial y\}$, where $y=0$ (respectively $x=0$) is any local analytic equation defining $E_i$ (resp. $E_j$).  
 \begin{proof} If $x'$ (resp. $y'$) are other equations, write $x'=ux, y'=vy,$ where $u,v$ are units in the local ring at $P_{ij}$; then, compare $x'\partial/\partial x'$ and $y'\partial/\partial y'$ with the previous choices, and reduce the coefficients modulo the maximal ideal.
 \end{proof}
 \end{lemma}

From now on, we restrict attention to 
graphs which are \emph{not} one of the following types:
\begin{enumerate}
\item A chain of smooth rational curves (= cyclic quotient singularity)
\item A cycle of smooth rational curves (= cusp singularity)
\item A smooth elliptic curve (= simple elliptic singularity)
\end{enumerate}
In every other case, $E$ contains at least one ``star'' curve $E_0$, with $2g_0+t_0>2.$  For such a curve, $h^0(\Theta_{E_0}(-t_0))=0$, hence $h^0(S\otimes\mathcal O_{E_0})=1.$  We recall a useful result:

\begin{proposition}(\cite{w})\label{w}  Assume that $H^0(S)\rightarrow H^0(S\otimes \mathcal O_{E_0})$ is surjective, where $E_0$ is a star curve.  Then $(V,0)$ is weighted homogeneous.
\begin{proof}  While this result is not stated explicitly in \cite{w}, a complete proof is found there from $(3.12)$ through $(3.16)$. One lifts a non-$0$ element of $H^0(S\otimes \mathcal O_{E_0})$ to a $D\in H^0(S)$, a derivation of the local ring of $V$.  A local argument along $E_0$ shows it is a non-nilpotent derivation, whence by a theorem of Scheja-Wiebe one has quasi-homogeneity.
\end{proof}
\end{proposition}

Let $R$ be the union of the rational end curves (with $g_i=0, t_i=1$), and $E'=E-R$ the union of the other curves.  $E'$ is connected and is the union of the stars plus rational curves with $t=2$.

\begin{lemma}\label{l1}  Suppose $F$ is a connected and reduced cycle in $E'$ containing a star $E_0$.  Then one has an inclusion into a one-dimensional space: $$H^0(S\otimes \mathcal O_F)\subset H^0(S\otimes \mathcal O_{E_0}) .$$
\begin{proof} Induct on the number of components of a connected $F'$ between $E_0$ and $F$.  Given such an $F'<F$, choose an $E_i$ in $F-F'$ which intersects $F'$.  Then $E_i\cdot F'>0$, and an easy check shows one has vanishing of $H^0$ of the first term in 
$$0\rightarrow S\otimes \mathcal O_{E_i}(-F')\rightarrow S\otimes \mathcal O_{F'+E_i}\rightarrow S\otimes \mathcal O_{F'}\rightarrow 0.$$
\end{proof}
\end{lemma}

A graph is called \emph{star-shaped} if it consists of one star out of which emanate strings of rational curves; it is the graph of a weighted homogeneous singularity.

\begin{lemma}\label{l4}  Suppose the graph of $X$ is star-shaped.  Then $h^0(S\otimes \mathcal O_{E'})=1.$ 
\begin{proof}
 As in \ref{l1}, the relevant connected cycles $F$ are constructed from the central star $E_0$ by adding smooth rational curves $E_i$ with $t_i=2$ with $F'\cdot E_i=1$.  Then  $$S\otimes \mathcal O_{E_i}(-F')=\mathcal O_{E_i}(-1)\oplus \mathcal O_{E_i}(-1),$$ so $h^0(S\otimes \mathcal O_{F'+E_i})=h^0(S\otimes \mathcal O_{F'}).$
 \end{proof}
\end{lemma}
\begin{lemma}\label{l5} Suppose the graph of $X$ is not star-shaped.  Then $h^0(S\otimes \mathcal O_{E'})=0.$
\begin{proof} Since the graph is not star-shaped, one concludes that either
\begin{enumerate}
\item $E'$ contains two stars $E_0$ and $E_0'$ connected by a (possibly empty) chain of rational curves with $t_i=2$
\item $E'$ contains one star $E_0$ and a chain of rational curves with $t_i=2$ and  intersecting $E_0$ at least twice.
\end{enumerate}
In either case, we claim that the union $F$ of the star and the rational curves in the chain satisfies $h^0(S\otimes \mathcal O_F)=0.$  The assertion then follows using the method of proof of Lemma \ref{l1}.

The result is clearest in the case where two stars $E_0$ and $E_0'$ meet at a point $P$.  In an appropriate affine open neighborhood $U$ of $P$, choose local functions $x$ and $y$ whose vanishing defines the two curves $E'_0$ and $E_0$, respectively, and so that $dx$ and $dy$ form a basis for the local $1$-forms.  Then $S$ is locally a free $\mathcal O(U)$-module with basis $x\partial/\partial x$ and $y\partial/\partial y$. The one-dimensional spaces of global sections of $S\otimes \mathcal O_{E_0}$ and $S\otimes \mathcal O_{E_0'}$ are of the form 
$ay\partial/\partial y$ and $bx\partial/\partial x$, respectively, where $a,b$ are constants.  These agree at $P$ and hence extend to a global section of $S\otimes \mathcal O_F$ only when $a=b=0.$

In the general case, we need to choose compatible ``coordinates'' on the components of $F$, as done in \cite{vt} (itself modeled closely on \cite{laufert}).  Denote the rational curves in the chain in order by $E_1,\cdots, E_r$, and let $E_{r+1}=E_0'$ be the star at the end.  Let $P_i=E_{i}\cap E_{i+1}$, $0\leq i\leq r$.  For $1\leq i\leq r$, the embedding of $E_i$ in $X$ is locally analytically the embedding in the normal line bundle of degree $-d_i$.  We choose coordinates so that $P_{i-1}$ corresponds to $\{0\}$ and $P_i$ to $\{\infty\}$.  Cover the scheme $2E_i$ by $2$ affines $$U_{i,1}=\text{Spec}\ k[x_i,y_i]/y_i^2$$
$$U_{i,2}=\text{Spec}\ k[x_i',y_i']/y_i'^2,$$
on whose intersection one has $$x_i'=1/x_i,\ \ \ y_i'=x_i^{d_i}y_i.$$ We may also assume that $y_i$ (and $y_i'$) are local equations for $E_i \subset 2E_i$, and that (possibly replacing $x_i$ by $x_i+y_ig(x_i)$, and similarly for $x_i'$) we may assume that the divisor of the predecessor $E_{i-1}\cap 2E_i\subset 2E_i$ has local equation given by $x_i=0$ (and similarly $x_i'=0$ at $\infty$).  In particular, we can assume that at each intersection point, the functions $x_i,y_i$ are restrictions (modulo higher order terms ) of local equations for the intersecting curves.  Furthermore, starting at $P_1$ and adjusting constants, we can assume that in the tangent space of $P_i$, we have $x_i'=y_{i+1}, \ y_i'=x_{i+1}$ (for $i<r$).  The standard exact sequence on $E_i$ $$0\rightarrow \mathcal O_{E_i}\rightarrow S\otimes \mathcal O_{E_i}\rightarrow \Theta_{E_i}(-t_i)\rightarrow 0$$  may be expressed on $U_{i,1}$  (since $t_i=2$) as $$0\rightarrow \{y_i\partial/\partial y_i\}\rightarrow\{y_i\partial/\partial y_i,\ x_i\partial/\partial x_i\}\rightarrow \{x_i\partial/\partial x_i\}\rightarrow 0,$$ and similarly on $U_{i,2}$.  The patching condition is
 $$x_i\partial/\partial x_i=-x_i'\partial/\partial x_i' +d_iy_i'\partial/\partial y_i'$$
 $$y_i\partial/\partial y_i=y_i'\partial/\partial y_i'.$$  For $1\leq i \leq r$, a global section of $S\otimes \mathcal O_{E_i}$ is of the form $$A_i y_i\partial/\partial y_i\ +B_i x_i\partial/\partial x_i=(A_i+d_iB_i)y_i'\partial/\partial y_i'\ -B_ix_i'\partial/\partial x_i',$$ for some constants $A_i,B_i.$

 For $1\leq i\leq r-1$, the two-dimensional space $S\otimes \C_{P_i}$ has the natural ordered basis  $$x_i'\partial/\partial x_i'=y_{i+1}\partial/\partial y_{i+1},\ y_i'\partial/\partial y_i'=x_{i+1}\partial/\partial x_{i+1}.$$  One similarly has an ordered basis at both $P_0$ and $P_r$.   Via patching, we have that
$$H^0(S\otimes \mathcal O_F)=\text{Ker}\ (\bigoplus_{i=0}^{r+1} H^0(S\otimes \mathcal O_{E_i})\rightarrow \bigoplus_{i=0}^{r} H^0(S\otimes \C_{P_i})).$$  Compatibility of the global sections above of the $H^0(S\otimes \mathcal O_{E_i})$ at $P_1,...,P_{r-1}$ means $$-B_i=A_{i+1},\ \ \ A_i+dB_i=B_{i+1},\ \  i=1,...,r-1.$$  A global section of the one-dimensional space $H^0(S\otimes \mathcal O_{E_0})$ is of the form $Bx\partial/\partial x$, where $x$ is a local equation of $E_0$ near $P_0$.  Therefore, its image in the space $S\otimes \mathbb C_{P_0}$ is $Bx_1\partial/\partial x_1$.  A section of $S\otimes \mathcal O_{E_1}$ patches compatibly if $A_1=0.$  Similarly, a global section of $H^0(S\otimes \mathcal O_{E_r})$ patches compatibly with a section of $H^0(S\otimes \mathcal O_{E_0'})$ if $A_r+d_rB_r=0.$ These $2r$ equations in the $A_i,B_i$ become $r$ equations in the $B_i$, namely 
$$-d_1B_1+B_2=0$$
$$B_1-d_2B_2+B_3=0$$
$$............$$
$$B_{r-1}-d_rB_r=0.$$

One recognizes the matrix of these equations as the intersection matrix of the cyclic quotient singularity whose resolution dual graph is that of the $r$ curves between $E_0$ and $E_0'$.  In particular, the determinant is $\pm n$, where one has an $n/q$ cyclic quotient.  Thus, the $B_i$ are all $0$.

Note that the same proof applies in case $E_0=E_0'$, except that one then has an additional condition that $B_1=-B_r$.
\end{proof}
\end{lemma}
\begin{proposition}\label{p6} Consider the minimal good resolution $(X,E)\rightarrow (V,0)$ of a normal surface singularity, with graph $\Gamma$, excluding the $3$ cases above.  Let $E'=E-R$, where $R$ is the union of the rational end curves, $r$ in number.  Then
\begin{enumerate}
\item if $\Gamma$ is star shaped, then $$h^1(S\otimes \mathcal O_{E'})=h^1(S\otimes \mathcal O_{E})=r+4h^1(\mathcal O_E)-3.$$ 
\item if $\Gamma$ is not star-shaped, then $$h^1(S\otimes \mathcal O_{E'})=h^1(S\otimes \mathcal O_{E})=r+4h^1(\mathcal O_E)-4.$$

\end{enumerate}

\begin{proof}  We claim that $$\chi(S\otimes \mathcal O_{E'})=4-4h^1(\mathcal O_E)-r,$$ where $r=\# R$ is the number of rational end curves.  For, Riemann-Roch for a rank $2$ vector bundle $G$ on $X$ restricted to a cycle $Z$ supported on $E$ states
$$\chi(G\otimes \mathcal O_Z)=-Z\cdot(Z+K)+Z\cdot \text{det}\ G.$$ (This may be easily deduced from the formula when $G$ is a line bundle.)  We let $Z=E'=E-R$ and $G=S$, and note that det($S)=-(K+E)$ and $E\cdot (E+K)=2h^1(\mathcal O_E)-2$. Now a small calculation establishes the claim.  We conclude $$h^1(S\otimes \mathcal O_{E'})=r+4h^1(\mathcal O_E)-4+h^0(S\otimes \mathcal O_{E'}).$$
Now apply Lemmas 2.4 and 2.5.

Finally, the short exact sequence $$0\rightarrow  \mathcal O_R(-E')\rightarrow \mathcal O_E\rightarrow  \mathcal O_{E'}\rightarrow 0$$
has as first term the direct sum of $\mathcal O(-1)$ for the rational end curves.  Tensoring with $S$, $H^1$ of the first term is $0$, whence the equality of $H^1$ of $S\otimes \mathcal O_E$ and $S\otimes \mathcal O_{E'}$. 
\end{proof}
\end{proposition}
\begin{remark}  The last short exact sequence also implies that in the star-shaped case, $H^0(S\otimes \mathcal O_E)\rightarrow H^0(S\otimes  \mathcal O_{E'})$ is a surjection onto a one-dimensional space, hence either space surjects onto $H^0(S\otimes  \mathcal O_{E_0})$, where $E_0$ is the central curve.
\end{remark}
\begin{remark}  In case $E$ consists of one smooth curve of genus $g>1$, one understands $h^1(S\otimes \mathcal O_E)=4g-3$ as corresponding to $3g-3$ deformations of the curve and $g$ deformations of the conormal bundle.
\end{remark}

We conclude with the useful
\begin{proposition}  Notation as above, \begin{enumerate}
\item $h^1(S)=h^1(S\otimes \mathcal O_E)+h^1(S(-E'))$ if $\Gamma$ is not star-shaped or $(V,0)$ is weighted-homogeneous
\item $h^1(S)=h^1(S\otimes \mathcal O_E)+h^1(S(-E'))-1$ if $\Gamma$ is star-shaped but $(V,0)$ is not weighted homogeneous.
\end{enumerate}
\begin{proof}  Using Proposition \ref{w} and the previous lemmas, one concludes that $H^0(S)\rightarrow H^0(S\otimes  \mathcal O_{E'})$ is the zero-map except when $(V,0)$ is not weighted homogeneous.  
\end{proof}
\end{proposition}

It is much easier to determine the cohomology of the restriction to $E$ of the determinant of $S$, namely $-(K_X+E)$.

\begin{proposition}\label{p2.10} Consider the minimal good resolution $(X,E)\rightarrow (V,0)$ of a normal surface singularity, with graph $\Gamma$, excluding the $3$ cases above.  Then $$h^1(-(K_X+E)\otimes \mathcal O_E)=r+3h^1(\mathcal O_E)-3.$$
\begin{proof} Since $-(K_X+E)\cdot E_i=-(2g_i-2+t_i)$ is \ $\leq 0$ except at rational end curves, and is strictly negative for at least one $E_i$, one easily concludes that $h^0(-(K_X+E)\otimes \mathcal O_{E'})=0.$  As at the end of the proof of the last Proposition, we have  $$h^0(-(K_X+E)\otimes\mathcal O_E)=h^0(-(K_X+E)\otimes \mathcal O_R(-E')),$$ which equals $r$.  The assertion now follows from Riemann-Roch, as $$\chi(-(K_X+E)\otimes \mathcal O_E)=(-3/2)E\cdot (E+K_X)=3(1-h^1(\mathcal O_E)).$$
\end{proof}
\end{proposition}

Recall that the Main Conjecture of \cite{mu} is an inequality about $$h^1(X,\mathcal O_X)-h^1(X,S_X)+h^1(X,-(K_X+E)).$$  For the corresponding expression for the cohomology of the restriction of these bundles to $E$, the previous results give a precise formula.
\begin{proposition}\label{p11} Let $(X,E)\rightarrow (V,0)$ be the minimal good resolution of a normal surface singularity, excluding the $3$ cases above.   Then $$h^1(\mathcal O_E)-h^1(S\otimes \mathcal O_E)+h^1(-(K_X+E)\otimes \mathcal O_E)=1-\delta,$$
where $\delta$ is $1$ if the graph is star-shaped, $0$ otherwise.
\end{proposition}
Thus, the Main Conjecture can be verified by explicit calculation in those cases for which $h^1$ of the bundles on $X$ agrees with $h^1$ of their restriction to $E$.  Since all one needs is an inequality, we have
\begin{corollary}\label{c1} Let $(X,E)\rightarrow (V,0)$ be the minimal good resolution of a normal surface singularity, excluding the $3$ cases above.  If $H^1(S(-E'))=0,$  then the Main Conjecture holds.
\begin{proof}  By Proposition \ref{p6},  $h^1(S)=h^1(S\otimes \mathcal O_{E'})=h^1(S\otimes \mathcal O_E).$  By the last Proposition, $$h^1(S)=h^1(\mathcal O_E)+h^1((-(K_X+E)\otimes \mathcal O_E)-1+\delta \leq h^1(\mathcal O_X)+h^1(-(K_X+E))-1+\delta,$$ so that $$h^1(\mathcal O_X)-h^1(S)+h^1(-(K_X+E))\geq 1-\delta \geq 0.$$
The Main Conjecture asserts that in the non-Gorenstein case, the left hand side is non-negative, and equals $0$ if and only the singularity is weighted homogeneous.  The inequality is now clear.  If the expression equals $0$, then $\delta=1$, so the graph is star-shaped.  For $E_0$ the star, consider the maps $$H^0(S)\rightarrow H^0(S\otimes \mathcal O_{E'})\rightarrow H^0(S\otimes \mathcal O_{E_0}).$$  The first map is surjective via the vanishing cohomology assumption, while the second is an isomorphism via Lemmas 2.3 and 2.4.  Proposition 2.2 then yields quasi-homogeneity.  Finally, it is proved in (\cite{mu}, Theorem $3.3$) that for a non-Gorenstein quasihomogeneous singularity, the expression is $0$.
\end{proof}
\end{corollary}
\begin{remark}  It is easy to see that $H^1(S(-E))=0$ implies $H^1(S(-E'))=0$, though the converse is false in general, even for rational singularities.
\end{remark}

\section{Vanishing theorems}
Let  $(X,E)\rightarrow (V,0)$ be a good resolution of a normal surface singularity, $\mathcal F$ a vector bundle on $X$.  We recall ``easy vanishing theorems'' (as in \cite{vt})  for  the local cohomology $$H^1_E(\mathcal F)=\lim_{\rightarrow}\ H^0(\mathcal F\otimes\mathcal O_Z(Z)).$$
This can be accomplished by a ``downward induction" on $Z$; if $E_i$ is in the support of $Z$, use the exact sequence $$0\rightarrow \mathcal F\otimes \mathcal O_{Z-{E_i}}(Z-E_i)\rightarrow \mathcal F\otimes \mathcal O_Z(Z)\rightarrow \mathcal F\otimes O_{E_i}(Z)\rightarrow 0.$$
\begin{proposition} (Vanishing Theorem)\label{vt1}  Suppose that for every $i$, $Z \cdot E_i<0$ implies $H^0(\mathcal F\otimes \mathcal O_{E_i}(Z))=0.$  Then $H^1_E(\mathcal F)=0.$

\end{proposition}
\begin{corollary}  If $L$ is a line bundle on $X$ with $L\cdot E_i \leq 0$, all $i$, then $H^1_E(L)=0.$

\end{corollary}

There is a slight refinement of the last Corollary which is occasionally useful.

\begin{proposition}\label{p2} If $L$ is a line bundle with $L\cdot E_i\leq 0$, all $i$, then either
\begin{enumerate}
\item $H^1_E(L(-E))=0$, or
\item $L\otimes \mathcal O_E\simeq \mathcal O_E$ and dim $H^1_E(L(-E))=1$.
\end{enumerate}
\begin{proof} One has the standard exact sequence $$0=H^0_E(L)\rightarrow H^0(L\otimes \mathcal O_E)\rightarrow H^1_E(L(-E))\rightarrow H^1_E(L)=0.$$
Suppose that $H^0(L\otimes \mathcal O_{E_i})=0$, some $i$.  By induction, we show that for a connected $F$ between $E_i$ and $E$, one has $H^0(L\otimes \mathcal O_F)=0.$  For, given such an $F$, choose an $E_j$ in $E-F$ which intersects $F$, and consider $$0\rightarrow H^0(L\otimes \mathcal O_{E_j}(-F))\rightarrow H^0(L\otimes \mathcal O_{F+E_j})\rightarrow H^0(L\otimes \mathcal O_F)\rightarrow 0.$$  By assumption, the degree of $L(-F)\otimes \mathcal O_{E_j}$ is negative, whence the assertion. 

 The only other possibility is that $L\otimes \mathcal O_{E_i}\simeq \mathcal O_{E_i}$, all $i$.  The same argument as above shows dim $H^0(L\otimes  \mathcal O_E)=1.$  This global section of $L\otimes \mathcal O_E$ induces an isomorphism on each $E_i$, so is an isomorphism itself.

\end{proof}
\end{proposition}
If $(V,0)$ has a rational singularity, then $Z\cdot (Z+K)<0,$ so every $Z$ contains an $E_i$ with $(Z+K)\cdot E_i <0$, i.e. $Z\cdot E_i <2-d_i.$  This yields the following refinements.
\begin{proposition}\label{p1} Assume $(V,0)$ is a rational singularity.  Suppose that for every $i$, $Z\cdot E_i<2-d_i$ implies $H^0(\mathcal F\otimes \mathcal O_{E_i}(Z))=0.$  Then $H^1_E(\mathcal F)=0.$
\end{proposition}
\begin{corollary}  Assume $(V,0)$ is a rational singularity and $L$ is a line bundle on $X$.  If $L\cdot E_i\leq d_i-2$, all $i$, then $H^1_E(L)=0.$
\end{corollary}
 
 We wish to study $H^1(S)$ and $H^1(-(K_X+E))$ in both the general and rational cases.  By duality, $h^1(\mathcal F)=h^1_E(\mathcal F^*\otimes K_X)$, so one can use the easy vanishing theorems, using the duals of the short exact sequences $$0\rightarrow \mathcal O_{E_i}\rightarrow S\otimes \mathcal O_{E_i}\rightarrow \Theta_{E_i}(-t_i)\rightarrow 0.$$ 

\begin{proposition}\label{gvt}  Let $(X,E)\rightarrow (V,0)$ be the minimal good resolution of a normal surface singularity.  Then
\begin{enumerate}
\item $H^1(S(-D))=0$ if for all $i$, $$D\cdot E_i \leq \text{min}\ \{2(2-2g_i)-t_i-d_i,\ 2-d_i\}.$$
\item $H^1(-(K_X+E)(-D))=0$ if for all $i$, $$D\cdot E_i\leq 2(2-2g_i)-t_i-d_i.$$
\end{enumerate} 
\begin{proof} For the first statement, note $h^1(S(-D))=h^1_E(S^*(K_X+D)).$  From the standard sequence for $S\otimes \mathcal O_{E_i}$, one has
 $$0\rightarrow K_{E_i}(t_i)\rightarrow S^*\otimes \mathcal O_{E_i}\rightarrow \mathcal O_{E_i}\rightarrow 0,$$

hence for any effective cycle $Z$, 
$$0\rightarrow K_{E_i}^{\otimes 2}(t_i
+d_i)(D+Z)\rightarrow S^*(K_X+D)\otimes \mathcal O_{E_i}(Z)\rightarrow K_{E_i}(d_i)(D+Z)\rightarrow 0,$$
Via Proposition $3.1$, one has only to check needed inequalities for $D\cdot E_i$ to guarantee that the first and third line bundles have negative degree whenever $Z\cdot E_i<0$.  One must handle separately the case of a rational end-curve ($g_i=0, t_i=1$), for which the second bound in the minimum is needed.

The second statement is proved similarly.

\end{proof}
\end{proposition}
\begin{proposition}\label{van}  Let $(X,E)\rightarrow (V,0)$ be the minimal good resolution of a rational surface singularity.  Then
\begin{enumerate}
\item $H^1(S(-D))=0$ if for all $i$, $$D\cdot E_i \leq \text{min}\ \{2-t_i,0\}.$$
\item $H^1(-(K_X+E)(-D))=0$ if for all $i$, $$D\cdot E_i\leq 2-t_i.$$
\end{enumerate} 
\begin{proof}  Same as the preceding proposition, but now using Proposition \ref{p1}.
\end{proof}
\end{proposition}

As indicated by Proposition \ref{p2} above, if one of the vanishing theorems holds for $D=nE$, one can sometimes conclude vanishing for $D=(n-1)E$, with only a mildly more restrictive condition.  
\begin{theorem} \label{vtt}  Let $(X,E)\rightarrow (V,0)$ be the MGR of a normal surface singularity.  Suppose that for all $i$, one has $$(*)\  d_i\geq 2(2g_i-2)+3t_i,$$ with strict inequality for at least one $i$.
Then 
\begin{enumerate}
\item $h^1(S(-E))=h^1(-(K_X+E)(-E))=h^1(\mathcal O_X(-E))=0.$
\item If $\Gamma$ is star-shaped, then $(V,0)$ is weighted homogeneous, and 
$h^1(S)=r+4h^1(\mathcal O_E)-3.$
\item If $\Gamma$ is not star-shaped, then $h^1(S)=r+4h^1(\mathcal O_E)-3.$
\item $h^1(\mathcal O_X)-h^1(S_X)+h^1(-(K_X+E))=1-\delta,$
where $\delta$ is $1$ if $(V,0)$ is weighted homogeneous, $0$ otherwise. 
\end{enumerate}
\begin{proof}   By hypothesis, $(2K_X+3E)\cdot E_i \leq 0,$ all $i$.  By Proposition \ref{p2}, this implies that $h^1_E(2K_X+2E)=0$, as long as the inequality is strict for some $i$.  Dually, $h^1(-(K_X+E)(-E))=0.$  Similarly, one can show $h^1(\mathcal O_X(-E))=0$ as long as $d_i\geq 2g_i-2+2t_i$, all $i$, with at least one strict inequality.

Applying Proposition \ref{gvt}(1) with $D=2E$, one has $h^1(S(-2E))=0$ as long as one has the inequalities $(*)$ and $d_i\geq 2t_i-2$.   This second inequality is a consequence of $(*)$ unless $g_i=0, t_i=1$; but in that case the extra inequality is $d_i\geq 0$, which is automatic.  Thus $h^1(S(-E))=h^1(S(-E)\otimes \mathcal O_E)$.  The dualizing sheaf of $E$ is the restriction of $K_X+E$, so by duality on $E$, the last space has dimension $h^0(S^*(K_X+2E)\otimes \mathcal O_E).$  For each $i$ we have the exact sequence 

$$0\rightarrow K_{E_i}^{\otimes 2}(3t_i
-d_i)\rightarrow S^*(K_X+2E)\otimes \mathcal O_{E_i}\rightarrow K_{E_i}(2t_i-d_i)\rightarrow 0.$$ The two line bundles are easily checked to have degree $\leq 0$ given $(*)$.  If for some $i$ we have $d_i>2(2g_i-2)+3t_i,$ both line bundles have strictly negative degree (one deals separately with the case $t_i=1, g_i=0$).  Therefore, $h^0(S^*(K_X+2E)\otimes \mathcal O_{E_i})=0.$ One can now apply the same induction trick on connected subcycles as in the proof of Proposition \ref{p2} to conclude that $h^0(S^*(K_X+2E)\otimes \mathcal O_{E})=0.$

By the cohomology vanishing of $(1)$, the expressions in $(2)$ and $(3)$ come from Proposition \ref{p6}; also, $(4)$ is given by Proposition \ref{p11}.  It remains to show that the graph is star-shaped if and only if the singularity is weighted homogeneous.  One direction is obvious; the other is proved as in Corollary \ref{c1}, replacing $E'$ there by $E$.
\end{proof}
\end{theorem}
\begin{remark}  The condition $(*)$ as well as the elimination of certain simple graphs means that none of the singularities in the Theorem could be Gorenstein.  The condition $(*)$ in case $E$ is one smooth curve is Grauert's well-known theorem that such a singularity is a cone.
\end{remark}
\begin{theorem}\label{rt} Let $(X,E)\rightarrow (V,0)$ be the MGR of a rational surface singularity.  Suppose that for all $i$, one has $$d_i\geq 2t_i-2.$$
\begin{enumerate}
\item  $h^1(S(-E))=h^1(-(K_X+E)(-E))=0.$
\item If $\Gamma$ is star-shaped, then $(V,0)$ is weighted homogeneous, and $h^1(S)=h^1(-(K_X+E))=r-3.$
\item If $\Gamma$ is not star-shaped, then $h^1(S)=h^1(-(K_X+E))-1=r-4.$
\end{enumerate}
\begin{proof}  $(1)$ follows directly from Proposition \ref{van} with $D=E$; the other statements follow as in the preceding proof.
\end{proof}
\end{theorem}
\begin{corollary}\label{rtc} For a rational singularity $(V,0)$ with $d_i\geq 2t_i-2$, all $i$, one has
\begin{enumerate}
  \item if $\Gamma$ is star-shaped, then dim $T^1_V=\sum_i(2d_i-3)+r-4.$
     \item if $\Gamma$ is not star-shaped, then dim $T^1_V=\sum_i(2d_i-3)+r-5.$
\end{enumerate}
\begin{proof}  That $d_i\geq 2t_i-2$ implies that $d_i\geq t_i$, hence the fundamental cycle is reduced.  Further, $d_i>t_i$ unless $d_i=t_i=2$.  Therefore, blowing up the maximal ideal yields only rational double point singularities.  According to the main formula of \cite{dej}, the dimension of $T^1$ minus the dimension of the Artin component  equals the multiplicity of the singularity minus $3$; there is no contribution from the other infinitely near points.  The formulas now follow.
\end{proof}
\end{corollary}

Of course, in general higher-order infinitesimal neighborhoods of $E$ can contribute to the three invariants in the Main Conjecture. 
 
  \begin{proposition}  Suppose the MGR $(X,E)\rightarrow (V,0)$ has $E$ a smooth curve of genus $g\geq 2$, whose conormal bundle $\mathcal O_E(-E)\equiv L$ has degree $d>2g-2$.  Then
  \begin{enumerate} 
  \item $h^1(\mathcal O_X)=g$
  \item $h^1(-(K_X+E))=3g-3+h^0(E,2K_E-L)$
  \item $h^1(S)\leq h^1(\mathcal O_X)+\ h^1(-(K_X+E))$, and equality is equivalent to quasi-homogeneity.
  \end{enumerate}
\begin{proof}  One computes directly that $h^1(S(-2E))=0$, using  that $$H^0(S^*(2E+K_X)\otimes \mathcal O_E(nE))=0, n\geq 1.$$   Similarly,  $h^1(-(K_X+E)(-2E))=h^1(\mathcal O_X(-E))=0.$  Thus, $H^0(S)\rightarrow H^0(S\otimes \mathcal O_{2E})$ is surjective, and $h^1(S)=h^1(S\otimes \mathcal O_{2E})$.  Now use the long exact sequence $$H^0(S\otimes \mathcal O_{2E})\rightarrow H^0(S\otimes \mathcal O_{E})\rightarrow H^1(S\otimes \mathcal O_E(-E))\rightarrow H^1(S\otimes \mathcal O_{2E})\rightarrow H^1(S\otimes \mathcal O_{E})\rightarrow 0.$$  The first map, into a one-dimensional space, is surjective if and only if $H^0(S)\rightarrow H^0(S\otimes \mathcal O_E)$ is surjective, which as mentioned earlier is equivalent to quasi-homogeneity.   One examines the same sequence with $S$ replaced by $-(K_X+E)$.  Everything now follows easily using Propositions \ref{p6} and \ref{p2.10}.
\end{proof}

\end{proposition}
One can show that all examples of this type are obtained by deforming the minimal resolution of the cone over $(E,L)$, by varying $E$ (in $3g-3$ ways) and $L$ (in $g$ ways); the cone has as well $h^1(E,\Theta \otimes L^{-1})$ non-conical equisingular deformations.

\section{Rational singularities}
    For rational singularities, there is an explicit topological formula for the second plurigenus $h^1(-(K_X+E))$.  
    
    \begin{theorem}\cite{mu},(4.4)  On the MGR of a rational surface singularity (not an RDP), let $Y$ be the smallest effective cycle satisfying $Y\cdot E_i\leq 2-d_i,$ all $i$.  If $Z=Y-E,$ then $$h^1(-(K_X+E))=Z\cdot (Z+3K)/2 +Z\cdot E.$$
  \end{theorem}  
      At a node, $E\cdot E_i=t_i-d_i>2-d_i$, so $Y$ has multiplicity at least 2 there; the same then applies for any neighbors with $t=2$.  Therefore, unless $Y=0$ (the RDP case) or $Y=E$ (cyclic quotient), one has $Y\geq E+E'.$ 
    
    \begin{proposition}\label{t'}  Exclude RDP's and cyclic quotients. Then $Y=E+E'$ iff $$(**)\  d_i\geq t_i+t'_i-2,$$ all $i$.  In this case, $h^1(-(K_X+E))=r-3,$ where $r$ is the number of ends of the graph.
    \begin{proof}  Examine $(E+E')\cdot E_i\leq 2-d_i$ and Theorem 4.1 (cf. also Proposition \ref{p6}).
    \end{proof}
    \end{proposition}

\begin{example}\label{4.3} The rational graph below satisfies  $(**)$ and has multiplicity $8$; as usual, the unmarked bullets are $-2's$.

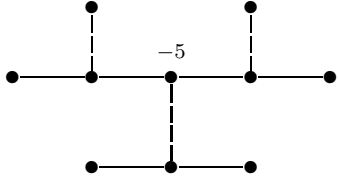
$$
\xymatrix@R=6pt@C=24pt@M=0pt@W=0pt@H=0pt{
\\
&&&\overtag{\bullet}{}{6pt}&
&\overtag{\bullet}{}{6pt}&&&\\
&&&\lineto[u]&&\lineto[u]&\\
&&&\lineto[u]&&\lineto[u]&\\
&&\undertag{\bullet}{}{2pt}\lineto[r]
&\undertag{\bullet}{}{2pt}\lineto[r]\lineto[u]
&\overtag{\bullet}{-5}{10pt}\lineto[r]&\undertag{\bullet}{}{2pt}\lineto[r]\lineto[u]
&\undertag{\bullet}{}{2pt}
&\\
&&&&\lineto[u]&\\
&&&&\lineto[u]&&\\
&&&&\lineto[u]&\\
&&
&\overtag{\bullet}{}{8pt}\lineto[r]
&\undertag{\bullet}{}{6pt}\lineto[r]\lineto[u]
&\undertag{\bullet}{}{6pt}&&&\\
&&&
&&&
&&&\\
}
$$ 
\medskip

\end{example}    

    Thus, for graphs satisfying $(**)$  (a generalization of the already considered case $d_i\geq2t_i-2$),  the Rational Conjecture states that $h^1(S)\leq r-3 $, with equality exactly in the quasihomogeneous case.  But  $h^1(S)\geq h^1(S\otimes\mathcal O_E)$, and the second term is $r-3$ or $r-4$, depending upon whether the graph is star-shaped or not.  So, in this case the Rational Conjecture is equivalent to the assertion $h^1(S)=h^1(S\otimes\mathcal O_E)=h^1(S\otimes\mathcal O_{E'})$.  In fact, there is a converse to Corollary \ref{c1}:
    
    \begin{proposition}\label{p4}  Let $(V,0)$ be a rational singularity whose graph satisfies $$\  (**) \ d_i\geq t_i+t'_i-2, \text{all}\  i.$$  Then the Rational Conjecture for $(V,0)$ is equivalent to $H^1(S(-E'))=0.$
    \begin{proof}  The first map in the exact sequence $$H^0(S)\rightarrow H^0(S\otimes\mathcal O_{E'})\rightarrow H^1(S(-E'))\rightarrow H^1(S)\rightarrow H^1(S\otimes \mathcal O_{E'})\rightarrow 0$$
    is the zero-map unless the singularity is weighted homogeneous, in which case it is surjective onto a one-dimensional space.
       The missed cases of RDP's and cyclic quotients are easily verified separately.
    \end{proof}
    \end{proposition}
    The next few sections will be devoted to proving the following result.
    
    \begin{theorem}\label{t1} Let $(V,0)$ be a rational singularity whose graph satisfies $$\  (**) \ d_i\geq t_i+t'_i-2, \text{all}\  i.$$ Then $H^1(S(-E')\otimes \mathcal O_E)=0.$
    \end{theorem}
 Unfortunately, at present we cannot conclude that $H^1(S(-E'))=0$ without a further hypothesis.   
    \begin{corollary}\label{c9} Suppose a rational singularity satisfying  $(**)$ has a reduced fundamental cycle, i.e. $d_i\geq t_i$, all $i$.  Then $H^1(S(-E'))=0,$ $h^1(S)=h^1(S\otimes \mathcal O_E)$,
 and the Rational Conjecture is true.
    \begin{proof}  $h^1(S(-(E+E'))=0$ by Proposition \ref{van}(1) and the hypotheses, as a simple calculation checks.  Thus, $h^1(S(-E'))=h^1(S(-E')\otimes \mathcal O_E)=0,$ by the Theorem.
    \end{proof}
    \end{corollary}
    
      Note Example \ref{4.3} is not covered by the Corollary; but see Example \ref{7.11} below.  Theorem \ref{t1} does not follow from an ``easy vanishing theorem,'' and is false in characteristic $p$ (Example \ref{pex}).

\section {Computing $H^1(S(-E')\otimes \mathcal O_E)$}
   There is an exact sequence $$0\rightarrow \mathcal O_E\rightarrow \oplus \mathcal O_{E_i}\rightarrow \oplus \mathbb C_{P_{ij}}\rightarrow 0,$$
 where the key map compares functions on adjacent curves $E_i$ and $E_j$ with their values at the intersection point $P_{ij}$.  (That is, on each $E_i$, at an intersection point $P_{ij}$ one sends a function $f$ to $\pm f(P_{ij})$, doing the opposite for $E_j$ at that point.) Tensoring with $S(-E')$ gives the important map 
 $$\Phi_E: \bigoplus_{E_i\subset E} H^0(S(-E')\otimes \mathcal O_{E_i})\rightarrow \bigoplus_{P_{ij}=E_i\cap E_j} H^0(S(-E')\otimes \mathbb C_{P_{ij}}).$$
 
 \begin{lemma}\label{Cok}\begin{enumerate}
 \item $H^1(S(-E')\otimes \mathcal O_{E_i})=0$ unless $d_i\leq t_i+t'_i-4.$
 \item If $d_i\geq  t_i+t'_i-3$ for all $i$, then $$\text{Coker}\  \Phi_E=
 H^1(S(-E')\otimes \mathcal O_E).$$
 \end{enumerate}
 \begin{proof} The first assertion is easily verified. For the second, tensor the short exact sequence with $S(-E')$ and take cohomology.
 \end{proof}
 \end{lemma}
The task is therefore to prove surjectivity of $\Phi_E$.  For $P_{ij}=E_i\cap E_j$, one gets contributions to the two-dimensional space $S(-E')\otimes \C_{P_{ij}}$ from $H^0(S(-E')\otimes \mathcal O_{E_i})$ and $H^0(S(-E')\otimes \mathcal O_{E_j})$.   We speak of the contribution of $E_i$ to its $t_i$ ``slots''.  The goal is to account systematically for contributions at the intersection points from the various edges.  We outline how to do this by an induction, using increasing sequences of connected subgraphs called \emph{cones}, handling one new intersection point at a time.

Return to graph language, and define a class of subtrees as follows: Let $v$ be a vertex of the graph $\Gamma$, with $p$ an adjacent edge. Define the \emph{cone} $\mathcal C(v,p)$ to be the connected component of $v$ in the graph $\Gamma-\{p\}$, plus the edge $p$ sticking out of it.  In other words, $\mathcal C(v,p)$ consists of all vertices and edges on the ``other side'' of $v$, away from $p$; but we keep the edge $p$ as well.  Thus, $\mathcal C(v,p)$ arises from adding $v$ and $p$ to $t_v-1$ other cones $\mathcal C(v_i,p_i)$, where the $p_i$ are the other edges emanating from $v$, with $v_i$ the other vertex of $p_i$; here $t_v$ is the valence of $v$.  We say that the $\{\mathcal C(v_i,p_i)\}$ are \emph{completed} by adding $v$ and $p$, forming $\mathcal C(v,p)$.
 
 In the Example  below,  $\mathcal C(v,p)$ consists of all vertices ``not below'' $v$, and is the completion of the three cones $\mathcal C(v_i,p_i)$, $i=1,2,3.$
$$
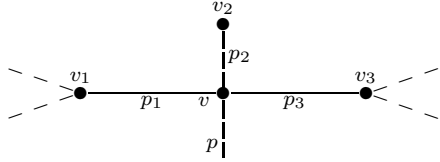
\xymatrix@R=6pt@C=24pt@M=0pt@W=0pt@H=0pt{
\\
&&&&\overtag{\bullet}{v_{2}}{6pt}
&&&&\\
&&&&\lineto[u]\righttag{}{p_2}{2pt}&&&\\
&&&&\lineto[u]&&&&\\
&&\overtag{\bullet}{v_{1}}{6pt}\dashto[ul]\dashto[dl]\lineto[r]&\overtag{}{p_1}{2pt}
\lineto[r]
&\lefttag{\bullet}{v}{2pt}\lineto[r]\lineto[u]\lineto[d]&\overtag{}{p_3}{2pt}\lineto[r]
&\overtag{\bullet}{v_3}{6pt}\dashto[dr]\dashto[ur]
&\\
&&&&\lineto[u]\lineto[d]&&&&&\\
&&&&\lefttag{}{p}{2pt}\lineto[u]\lineto[d]&&&&\\
&&&&\lineto[u]&&&&\\
}
$$ 
One may form sequences of cones as follows: in Round $1$, consider cones consisting of an end and its intersection point with its one neighbor.  In Round $2$, consider those $v$ all but one of whose neighbors are ends, with $p$ the edge leading to the other neighbor; these arise as the ends of $\Gamma$ minus the original ends. Form the corresponding cones $\mathcal C(v,p)$.  In Round $3$, again consider the ends of $\Gamma$ minus all the vertices in a previous cone; these will have the property that all but one of their neighbors are vertices in earlier cones.  Form a new cone by adding the missing edge $p$.  That is, we \emph{complete} previously considered cones.  At each Round, one completes at least one new cone.  Eventually, one is left with a graph with one node, all of whose neighbors occur in previously chosen cones.  We call this the \emph{terminal} situation.

Alternatively, one may start with any interior node $v$, whose edges $p_1,\cdots,p_t$ lead to vertices $v_1,\cdots,v_t$, and consider it the terminal stage of $t$ cones $\mathcal C(v_i,p_i)$.  Then, take apart each $\mathcal C(v_i,p_i)$ by reversing the completion process.  But to do our induction, one starts at the ends, and works ones way up to the terminal node $v$.

Each pair $(v,p)$ corresponds to a curve $E_0$ and an intersection point $P_0$, and the cone $\mathcal C(E_0,P_0)$ is the union of all curves leading away from $E_0$ via the $t_0-1$ intersection points other than $P_0$ (but with $P_0$ a distinguished point.)  From now on, we shall use curve (rather than node) notation.

We outline the induction process using increasing sequences of cones.
Suppose $E_0$ is an end-curve, intersecting at $P_0$ with another curve $E_1$.  We will show below that $$H^0(S(-E')\otimes  \mathcal O_{E_0})\rightarrow H^0(S(-E')\otimes \C_{P_0})$$
is an inclusion of a one-dimensional space; so in showing surjectivity of $\Phi_E$ we will have to account for the missing dimension at $P_0$ by using a contribution from $E_1$.  At the next step, change notation so that $E_0$ has valency $t$ and intersects $t-1$ end-curves at $P_1,\cdots,P_{t-1}$; if $P_0$ is the remaining intersection point, we are considering the cone $\mathcal C(E_0,P_0)$.  We need $H^0(S(-E')\otimes  \mathcal O_{E_0})$ to have enough sections to account for the missing dimensions at $P_1,\cdots,P_{t-1}$; that is, for each of these $P_i$ we need a section that vanishes at the other $t-2$ points and contributes the needed dimension at $S(-E')\otimes  \C_{P_i}$.  That would take care of the desired surjectivity at $t-1$ intersection points.  Optimally, we would also like the space of sections of $H^0(S(-E')\otimes  \mathcal O_{E_0})$ which vanish at these $t-1$ points to map onto $S(-E')\otimes C_{P_0}$.  Then we would not have to worry about these $t$ points for the rest of the induction.  We'll call this the Type I case.  But we will be satisfied if we can find a section vanishing on the $t-1$ points and giving a ``useful'' element of $S(-E')\otimes \C_{P_0}$; call this Type II.  The general case will consist of completing cones for which appropriate surjectivity results have been established.  At the last step, one has an $E_0$ all of whose $t$ intersection points arise from previously considered cones.

Recall that if $E_i$ is locally defined by $y=0$, $E_j$ by $x=0$, then $S\otimes \C_{P_{ij}}$ has a natural ordered basis $x\partial/\partial x, y\partial/\partial y$ which is independent of the choice of $x$ and $y$.  Multiplying by a local equation of $E'$ (either $x$, $y$, or $xy$) gives an ordered basis for $S(-E')\otimes \C_{P_{ij}}$, so that elements are given by an ordered pair of numbers $(a,b)$; this is equal (up to a scalar multiplication) to the element $(b,a)$ viewed from considering $E_j\cap E_i$.

Returning to cohomological considerations, for a cone $\mathcal C(E_0,P_0)$, we consider two natural ``evaluation'' maps:

$$\Phi_{E_0,P_0}:\bigoplus_{E_i\subset \mathcal C} H^0(S(-E')\otimes \mathcal O_{E_i})\rightarrow \bigoplus_{P_{ij}=E_i\cap E_j} H^0(S(-E')\otimes \mathbb C_{P_{ij}})$$
$$\Psi_{E_0,P_0}:\bigoplus_{E_i\subset \mathcal C} H^0(S(-E')\otimes \mathcal O_{E_i})\rightarrow \bigoplus_{P_{ij}=E_i\cap E_j} H^0(S(-E')\otimes \mathbb C_{P_{ij}})\oplus H^0(S(-E')\otimes \C_{P_0})$$
 
 In other words, for $\Phi$ we consider all points of intersection of curves in $\mathcal C$, $t_0-1$ of which are on $E_0$; for $\Psi$, we also consider evaluation at the additional point $P_0$.  Clearly, 
 $$\Phi_{E_0,P_0}=\pi \cdot \Psi_{E_0,P_0},$$ where $\pi$ is projection off the two-dimensional direct summand $H^0(S(-E')\otimes \C_{P_0})$.  We will study the cokernel of these maps for our judiciously chosen increasing sequence of cones, ultimately concluding the desired vanishing result.
 
 It will turn out that if $d_i\geq t_i+t'_i-2$ for all $i$, then all cones (except from end-curves) will have one of two properties:
 
 \begin{definition} $\mathcal C(E_0,P_0)$ is \emph{Type I} if $\Psi_{E_0,P_0}$ is surjective.
 \end{definition}
 \begin{definition} $\mathcal C(E_0,P_0)$ is \emph{Type II} if 
 \begin{enumerate}
 \item $\Psi_{E_0,P_0}$ has image of codimension $1$
 \item $\Phi_{E_0,P_0}$ is surjective
 \item $\text{Im}\ \Psi\ \cap\ \text{Ker}\ \pi$ contains an element in $S(-E')\otimes \C_{P_0}$ with coordinates $(1,-b)$, where $b\geq 1$.
 \end{enumerate}
 \end{definition}

 This will be done by a variant of the method in Section 2; local coordinates on $2E$ will be chosen in the same way, but now one allows that $t_i$ may be greater than $2$.

\section{{The cokernels of $\Phi$ and $\Psi$ for cones}}

In considering a cone $\mathcal C(E_0,P_0)$, we choose coordinates for the curves as in Section 2.  To simplify notation, we write $d=d_0$, $t=t_0$, and $t'=t'_0$. An $E_0$ will be defined in a first chart by $y=0$, with $P_0$ given by $x=0$, and the other intersection points (if any) $P_1,\cdots,P_{t-1}$ given by $x=a_1,\cdots,a_{t-1}$ (no intersection points at $\infty$).  Then $S\otimes \mathcal O_{E_0}$ will be generated by the images of $x\prod_{j=1}^{t-1} (x-a_j)\partial/\partial x$ and $y\partial/\partial y.$ The second chart, defined by $x',y'$, will be as in Section 2.  We write out the elements of $H^0(S(-E')\otimes \mathcal O_{E_0})$ and compute their images in the various ``slots'', i.e., the two-dimensional spaces $S(-E')\otimes \C_{P_j}$.  Note while $S\otimes \C_{P}$ has a natural ordered basis (Lemma 2.1), the ordered basis for $S(-E')\otimes \C_{P}$ is unique up to a scalar multiple.

Note that the bundles we consider satisfy
\begin{enumerate}
\item $t>1$: $S(-E')\otimes \mathcal O_{E_0}\equiv \mathcal O_{E_0}(d-t')\oplus \mathcal O_{E_0}(d-t-t'+2)$
\item $t=1$: $S(-E')\otimes \mathcal O_{E_0}\equiv \mathcal O_{E_0}(-1)\oplus \mathcal O_{E_0}$.
\end{enumerate}
In the ``easy case," $H^1(S(-E')(-\sum_{i=0}^{t-1}P_i)\otimes \mathcal O_{E_0})$ vanishes, hence $$H^0(S(-E')\otimes \mathcal O_{E_0})\rightarrow \bigoplus _{i=0}^{t-1} H^0(S(-E')\otimes \mathbb C_{P_{i}})$$ is surjective, guaranteeing that $E_0$ fills all the slot entries at its $t$ intersection points.  We start with

\begin{lemma} \label{s} For a curve $E_0$ with $t=2$, suppose $t'=1$ or $d\geq 3.$  Then
$$H^0(S(-E')\otimes \mathcal O_{E_0})\rightarrow \bigoplus_{i=0}^{1} H^0(S(-E')\otimes \C_{P_i}) $$ is surjective.
\end{lemma}
\begin{lemma} \label{s2} Consider a cone $\mathcal C(E_0,P_0)$ consisting of a string of curves starting with an end-curve, for which $t_0=2$.  Then the cone has Type I.
\begin{proof}  Starting from the end-curve, name the curves in the string as $E_r, \cdots, E_1, E_0$, with intersection points $P_i=E_i\cap E_{i-1}, i=1,\cdots,r$. Then by Lemma \ref{s}, $E_{r-1}$ fills up both of its slots, so $\mathcal C(E_{r-1},P_{r-1})$ has Type I.  Moving along the chain, a later $E_i$ need only fill the slot at $P_i$, which is automatic because $H^1(S(-E')\otimes \mathcal O_{E_i}(-P_i))=0$. \end{proof}
\end{lemma}
The general case requires more delicate argument.
\begin{lemma}\label{s3} Suppose $E_0$ is an end-curve, with self-intersection $-d$ and intersection point $P_0$.  Then $H^0(S(-E')\otimes \mathcal O_{E_0})$ is one-dimensional, with basis $x\cdot x\partial/\partial x - dx\cdot y\partial/\partial y$, whose image in $S(-E')\otimes \C_{P_0}$ is $(1,-d)$ in the ordered basis $x\cdot x \partial/\partial x, x\cdot y\partial/\partial y.$
\begin{proof} In the first chart, $E'$ is defined by $x=0$, so a section of $S(-E')\otimes \mathcal O_{E_0}$ is of the form $$A(x)x\cdot x\partial/\partial x + B(x)x\cdot y\partial/\partial y.$$
In the other chart, where $E'$ is empty, this becomes $$-A(1/x')\partial/\partial x' +\{(d/x')A(1/x') + 1/x'B(1/x')\} y'\partial/\partial y'.$$
To be a global section requires that $A$ is a constant, say $1$, in which case $B$ must be the constant $-d$.
\end{proof}
\end{lemma}
Thus, a Round $1$ cone $\mathcal C(E_0,P_0)$ could be considered of Type II, but we view these separately.

Now consider a cone $\mathcal C(E_0,P_0)$, where $P_0$ does not intersect an end-curve.  We write out the global sections of $S(-E')\otimes \mathcal O_{E_0}$.  Use coordinates for which $P_0$ is given by $x=0$, and the other intersection points $P_i$ are at $x=a_i,\ i=1,\cdots,t-1.$  Assume further that the last $t-t'$ of these are the points intersection with end-curves.  $S\otimes \mathcal O_{E_0}$ is generated on the first chart by $x\prod (x-a_i)\partial/\partial x$ and $y\partial/\partial y$, and on the second by $\prod(x'-(1/a_i))\partial/\partial x'$ and $y'\partial/\partial y'$. An equation for $E'$ is given by $z=y\cdot x\prod_{i=1}^{t'-1} (x-a_i)$ in the first chart, and $z'=y'\prod_{i=1}^{t'-1} (x'-(1/a_i))$ on the second.  For convenience, let $\alpha = \prod_{i=1}^{t'-1} (-a_i), \ \beta=\prod_{j=1}^{t-1}(-a_j)$.  A calculation yields the

\begin{lemma}\label{L} In the coordinates above, the sections of $S(-E')\otimes \mathcal O_{E_0}$ are written in the two charts as $$A(x)zx\prod_{i=1}^{t-1} (x-a_i)\partial/\partial x + B(x)zy\partial/\partial y=$$ 
$$z'\prod(x'-(1/a_i))\partial/\partial x'\{\alpha \beta(-x'^{d-t-t'+2}A(1/x'))\}+$$
$$z'y'\partial/\partial y'\{\alpha \beta \prod (x'-(1/a_i))dx'^{d-t-t'+1}A(1/x')+\alpha x'^{d-t'}B(1/x')\}.$$
These are global sections exactly when
\begin{enumerate}
\item $A(x)$  is a polynomial of degree $d-t-t'+2$, with coefficient of $x^{d-t-t'+2}$ denoted $C$
\item $B(x)$ is a polynomial of degree $d-t'+1$, with coefficient of $x^{d-t'+1}$ denoted $C'$
\item $C'=-dC.$
\end{enumerate}
\begin{proof} Convert to coordinates in the second chart, carefully.
\end{proof}
\end{lemma}

We record that for a global section as above, the element induced in the slot at $P_j$ for $j>0$ has coordinates $$(A(a_j)a_j\prod_{i\neq j}(a_j-a_i), B(a_j)),$$ while the coordinates at $P_0$ are $$(A(0)\beta, B(0)).$$ 

\begin{lemma}\label{s5}  Suppose $d=t=t'=2$, with $P_1$ given by $a_1=1$.  Then the general element of $H^0(S(-E')\otimes \mathcal O_{E_0})$ is of the form $Azx(x-1)\partial/\partial x +(B-dAx)zy\partial/\partial y,$ with $A,B$ arbitrary.  Its image in $S(-E')\otimes \C_{P_i}$ is $(-A,B)$ for $i=0$, $(A,B-dA)$ for $i=1$.
\end{lemma}
Excluding for the moment the case that $\Gamma$ is star-shaped with length one arms, one can form Round 2 cones.

\begin{lemma}\label{beg} Suppose $\mathcal C(E_0,P_0)$ is a Round $2$ cone; thus $t'=1$.  Then either
 \begin{enumerate}
\item $d\geq t$ and $\mathcal C$ is Type I, or
\item $d=t-1$ (so $t\geq 3$) and $\mathcal C$ is Type II.
\end{enumerate}
\begin{proof} Assume first that $d\geq t$.  Then in the notation of Lemma \ref{L}, we certainly have global sections with $A(x)=0$ and $B(x)$ a polynomial of degree $t-1$.  In particular, for every $P_i$, $i\geq 0$, we can choose $B(x)$ to be a polynomial taking value $1$ at $P_i$ and $0$ at $P_j, j\neq i$.  For $i\geq 1$, this contributes  $(0,1)$ in the coordinates of $S(-E')\otimes \C_{P_i}$, and $(0,0)$ in the $t-1$ other slots.  As the contribution from the corresponding end-curve in this slot is $(-d_i,1)$, we conclude that $\text{Im}\ \Psi_{E_0,P_0}$ contains the direct summand $\oplus_{i=1}^{t-1}S(-E')\otimes \C_{P_i}$.  On the other hand, we have $H^0(S(-E')\otimes \mathcal O_{E_0})\rightarrow S(-E')\otimes \C_{P_0}$ is surjective since $H^1(S(-E')\otimes \mathcal O_{E_0}(-P_0))=0.$  Thus, $\Psi_{E_0,P_0}$ is surjective, and the cone is of Type I. 

Next suppose $d=t-1$. There is a global section with $A(x)=0$ and $B(x)$ a polynomial of degree $t-2$.  Thus for every $1\leq i\leq t-1$, we can find a polynomial $B_i(x)$ which vanishes at all the $P_j$ except for $P_i$ and $P_0$. Combining with the contribution from the end-curves, we conclude that $\Phi_{E_0,P_0}$ is surjective.  To find an element in $\text{Im}\ \Psi \cap \text{Ker}\ \pi$, we consider the global section of $S(-E')\otimes \mathcal O_{E_0}$ with $A(x)=1$ and
$$B(x)=-d\prod_{i=1}^{t-1}(x-a_i)-\sum_{i=1}^{t-1}(a_i/d_i)\prod_{k\neq i}(x-a_k).$$
For $j>0, B(a_j)=-(a_j/d_j)\prod_{k \neq j}(a_j-a_k)$, so by the results above the contribution of the global section at this slot has coordinates $$(a_j\prod_{k\neq j}(a_j-a_k),\ -(a_j/d_j)\prod_{k \neq j}(a_j-a_k)).$$
This is a non-$0$ multiple of the section $(-d_j,1)$, which can be matched by a contribution from the corresponding end-curve.  Furthermore, $B(0)=\beta(-d+ \sum (1/d_i))$, hence the contribution of the section at $P_0$ has coordinates
$$(\beta,\beta(-d+\sum(1/d_i))).$$ Therefore, subtracting off contributions from the $t-1$ end-curves gives an element in the image of $\Psi$ whose only non-$0$ entries are at $P_0$, and it is a multiple of $(1,-(d-\sum(1/d_i))).$  But as $d_i\geq 2$, we conclude that $$d-\sum(1/d_i))\geq (t-1)-(t-1)/2=(t-1)/2  \geq 1$$ (as $t\geq 3)$.  Therefore, the corresponding cone is of Type II.

\end{proof}
\end{lemma}
We can now consider a cone $\mathcal C(E_0,P_0)$ in an arbitrary round, postponing the terminal situation for which the previously considered cones involve all but one curve.

\begin{lemma}\label{midround} Consider a cone $\mathcal C(E_0,P_0)$, with intersection points at $P_1,\cdots ,P_{t-1}$ coming from previously considered Rounds.  Assume as before that $d\geq t+t'-2.$  Then this cone has Type I unless $d=t+t'-2$ and $P_1,\cdots,P_{t-1}$ come from end-curves or Type II cones.  In that case, the cone has Type II.
\begin{proof} Recall that there is a global section of $S(-E')\otimes  \mathcal O_{E_0}$ with $A(x)=0$ and $B(x)$ a polynomial of degree $d-t'$.  Suppose that some number $t^*$ of the $t-1$ points come from Type II cones or end-curves.  Then for any of these points, there is a polynomial $B_i(x)$ of degree $t^*$ which is non-$0$ at that point but vanishes at the other $t^*-1$ other points and also at $P_0$.  If $d-t'\geq t^*$, such polynomials can be used to construct global sections, and so we can conclude surjectivity of $\Psi$ exactly as in the proof of the previous Lemma.  In other words, the cone is of Type I unless $d-t' <t^*$.  But by hypothesis $d-t'\geq t-2$; so the only case not covered is that $t^*=t-1$ and $d=t+t'-2$.

Considering this remaining case, one can as in the above Lemma choose for every $P_i \ (i>0)$ a polynomial of degree $d-t'=t-2$ which vanishes exactly at all $P_k, k\neq 0,i.$  As in the last Lemma, one concludes surjectivity of $\Phi$. It remains to produce an element in $\text{Im}\ \Psi \cap \text{Ker}\ \pi$ whose coordinates in the $P_0$ slot are of the form $(1,-u)$, where $u$ is a rational number $\geq 1$.  For this, we choose the global section as in the proof of the last Lemma, except that we must match the contribution $(1,-u_i)$ at $P_i$, where $u_i\geq 1$.  This choice will produce a contribution in the slot at $P_0$ of $(1,-d+\sum (1/u_i)).$  Thus, it remains only to verify that $$d-\sum_{i=1}^{t-1}(1/u_i)\geq 1.$$
As $u_i\geq 1$ and $d=t+t'-2$, the inequality follows easily unless $t'=1$.  But that means all $P_i \ (i>0)$ come from end-curves, and that case was handled in the preceding Lemma.
\end{proof}
\end{lemma}

We are ready to consider the terminal situation.  Suppose $E_0$ is a curve with intersection points $P_1,\cdots,P_t$ intersecting with other curves $E_1,\cdots,E_t$, and the corresponding cones $\mathcal C(E_i,P_i)$ have already been shown to be end-curves or of Type I or Type II as above.

\begin{lemma}  In the situation above, the map $\Phi_E$ is surjective, hence $$H^1(S(-E')\otimes \mathcal O_E)=0.$$
\begin{proof} We need to produce global sections of $S(-E')\otimes \mathcal O_{E_0}$ which account for missing entries in slots coming from Type II cones and end-curves.  Assume there are $t^*\leq t$ of these.  As above, if $d-t'\geq t^*-1$ we can find at each of these points a suitable global section with slot entry $(0,1)$ there, but vanishing at the other $t^*-1$ points.  This suffices to prove the surjectivity of $\Phi_E$ in that case.

But $d-t' < t^*-1$ only when $t^*=t$ and $d=t+t'-2$.  However, we claim that this cannot happen, because the graph is that of a rational singularity.  Consider the cycle $Z=E+E'$.  For an end-curve $E_i$, we have $Z\cdot E_i=2-d_i$, hence $(Z+K)\cdot E_i=0.$  The condition $t^*=t$ means that every cone along the way has been of Type II; for any other curve $E_i$ not $E_0$, we have $d_i=t_i+t_i'-2.$  As $Z\cdot E_i=t_i+t'_i-2d_i$, we conclude that also in this case $(Z+K)\cdot E_i=0.$  By rationality, $Z\cdot (Z+K)\leq -2$, so we must have $(Z+K)\cdot E_0 <0$.  But this says $t+t'-d-2<0$, contradicting the hypothesis.
\end{proof}
\end{lemma}
\begin{remark}\label{r1}  Note that if the fundamental cycle is reduced, i.e., $d_i\geq t_i$ for all $i$, then there are no curves of Type II; one has only end-curves and Type I curves.  By Lemma \ref{beg}, this is clear at Round $2$.  In a later Round, by Lemma \ref{midround} the only new Type II case would occur if there were $t-1$ end-curves; but that case was handled in the previous Round.
\end{remark}

\section{Some sharpened results}
We have shown that if $d_i\geq t_i+t'_i-2$, all $i$, then $H^1(S(-E')\otimes \mathcal O_E)=0$.  As clear from Lemma 5.1, vanishing is not possible if some $d_i\leq t_i+t'_i-4.$
In this Section we discuss some vanishing for graphs with one or more vertices satisfying  $$d=t+t'-3.$$  As $d\geq t-1$, we have $t'\geq 2.$ According to Lemma 6.2, on such a curve a global section of $S(-E')\otimes \mathcal O_{E_0}$ has $A(x)=0$ and $B(x)$ a polynomial of degree $t-3.$  So, for any set of the $t-3$ intersection points, one may chose a section vanishing at all of them, and giving a non-zero contribution of the form $(0,\cdot)$ at each of the other $3$ points.  We easily conclude the following two results.
\begin{lemma} \label{e1} Consider a cone $\mathcal C(E_0,P_0)$ so that $E_0$ satisfies $d=t+t'-3$.  Suppose that at least $2$ of the vertices $P_1,\cdots,P_{t-1}$ correspond to Type I cones.  Then $\mathcal C(E_0,P_0)$ has Type II, except that the contribution at $P_0$ is $(0,1)$.
\end{lemma}
\begin{lemma}\label{e2}  At the terminal stage, suppose $E_0$ satisfies $d=t+t'-3$, and at least two of the $t$ intersection points come from Type I cones.  Then $\Phi_E$ is surjective, hence $H^1(S(-E')\otimes \mathcal O_E)=0$.
\end{lemma}  
\begin{corollary} For a star-shaped rational graph whose central curve satisfies $d\geq t+t'-3$, one has $H^1(S(-E')\otimes \mathcal O_{E})=0$.
\begin{proof} Since $d\geq t-1$, one has in the exceptional case $d=t+t'-3$ that $t'\geq 2$, so there are at least $2$ Type I strings.
\end{proof}
\end{corollary}
\begin{lemma}\label{e3}  Suppose the graph contains two curves $E_0$ and $E'_0$ satisfying $d=t+t'-3$, connected by a (possibly empty) string of rational curves, while all other curves satisfy $d\geq t+t'-2$.  For the $t_0-1$ intersection points of $E_0$ not pointing towards $E_0'$, assume as in Lemma \ref{e1} that at least two correspond to Type I cones; make the same assumption for $E'_0$.  Then $\Phi_E$ is surjective, hence $H^1(S(-E')\otimes \mathcal O_E)=0).$
\begin{proof}  By assumption, $E_0$ is connected at some $P_0$ by a string of curves with $t=t'=2$ to some $P_0'\in E_0'$. If $P_0=P'_0$, then  via Lemma \ref{e1} the contributions from the two curves to $S(-E')\otimes \C_{P_0}$ span the whole space, so $\Phi_E$ is surjective.  

So, suppose $E_0$ is joined to $E'_0$ by a chain of $r$ rational curves $E_1, \cdots,E_r$, with $P_i=E_i\cap E_{i+1}$ ($i<r$).   Assume first that all intermediary curves $E_i$ have $d_i=2$.  Moving from $E_0$ towards $E_0'$, we show that all intermediate cones $\mathcal C(E_i,P_i)$ are of Type II, except that the extra contribution at $P_i$ is $(0,1)$.   By Lemma \ref{s5}, the global sections of $S(-E')\otimes \mathcal O_{E_i}$ give a contribution of $(-A_i,B_i)$ at $P_{i-1}$, and $(A_i,B_i-2A_i)$ at $P_{i}$.  So, use $A_i=-1, B_i=0$ to fill in the slot at $P_{i-1}$, then use $A_i=0$ and $B_i=1$ to make the contribution $(0,1)$ at $P_i$.   At the last stage, $E'_0$, which already had a $(1,0)$ at $P'_0$, now receives a $(0,1)$ from the last curve in the string.

If some $E_i$ satisfies $d_i\geq 3$, then by Lemma \ref{s}  $H^0(S(-E')\otimes \mathcal O_{E_i})$ maps onto the sum of the two spaces $S(-E')\otimes \C_P$, $P$ an  intersection point of $E_i$.  By Lemma \ref{s5}, then each curve adjacent to $E_i$ maps onto the space $S(-E')\otimes \C_P$ for $P$ the outer point in the direction away from $E_i$.  Continuing in this way gives the desired surjectivity, without even using the contributions of $E_0$ and $E'_0$ at $P_0$ and $P'_0$.

\end{proof}
\begin{remark}\label{r1} In the previous Lemma, if one assumes that the fundamental cycle is reduced, then as already mentioned there are no Type II curves except end-curves.  But among the $t_0-1$ neighbors of $E_0$ are at least $2$ non-end-curves, as $d-t_0=t'_0-3\geq 0$.  So, the conditions on intersection points are automatically satisfied. 
\end{remark}
\end{lemma}
We illustrate this case with several examples.  The first is known by \cite{laufert} to be a taut singularity, which fact may be deduced in two steps; the first is:
\begin{example}\label{ex1}  For $a,b\geq 3$, a singularity with resolution graph 
$$
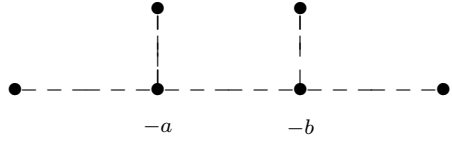
\xymatrix@R=12pt@C=24pt@M=0pt@W=0pt@H=0pt{
&&{\bullet}\dashto[dd]&&{\bullet}\dashto[dd]\\
&&\dashto[u]&&\dashto[u]\\
{\bullet}\dashto[r]&&
\undertag{\bullet}{-a}{3pt}\dashto[l]\dashto[u]\dashto[r]&\dashto[r]&
\undertag{\bullet}{-b}{3pt}\dashto[r]&\dashto[r]&{\bullet}\\
&&&&\\
}
$$
satisfies $H^1(S(-E')\otimes \mathcal O_E)=0.$
\begin{proof}Assume that all end-strings emanating from the two nodes $E_1$ and $E_2$ contain at least $2$ curves; the other cases are similar or easier.  Then Lemma \ref{e3} applies.
\end{proof}
\end{example}
\begin{remark} Note that in the above example, the plurigenus $h^1(-(K_X+E))$ may have dimension much bigger than $h^1(S)=0$.  For instance, if all non-nodal curves are $-2's$, and the outward end-strings each have length $n$, then the plurigenus equals $n$.
\end{remark}
\begin{example}  Singularities with the graph below have $H^1(S(-E')\otimes \mathcal O_E)\neq 0$, hence $h^1(S)>h^1(S\otimes \mathcal O_E))=1.$
$$
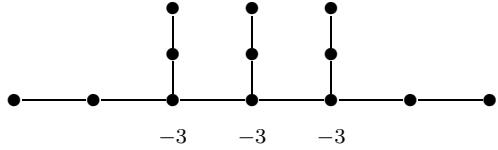
\xymatrix@R=12pt@C=24pt@M=0pt@W=0pt@H=0pt{
&&{\bullet}\lineto[d]&{\bullet}\lineto[d]&{\bullet}\lineto[d]\\
&&{\bullet}\lineto[d]&{\bullet}\lineto[d]&{\bullet}\lineto[d]\\
{\bullet}\lineto[r]&{\bullet}\lineto[r]&
\undertag{\bullet}{-3}{3pt}\lineto[l]\lineto[u]\lineto[r]&\undertag{\bullet}{-3}{3pt}\lineto[r]&
\undertag{\bullet}{-3}{3pt}\lineto[r]&{\bullet}\lineto[r]&{\bullet}\\
&&&&\\
}
$$
\begin{proof} The $3$ curves $E_i$ (corresponding to the $3$ nodes) each satisfy $d_i=t_i+t'_i-3$, so the sum of dimensions of the spaces of sections of $S(-E')\otimes \mathcal O_{E_i}$ equals $3$. But only they can contribute to the two two-dimensional spaces $S(-E')\otimes \C_{P}$ at the two intersection points (the edges joining the nodes).  So $\Phi_E$ cannot be surjective.  

This example is still consistent with the Rational Conjecture, as the plurigenus equals $4$.
\end{proof}
\end{example}
\begin{theorem}\label{t3}  Suppose a rational singularity, with reduced fundamental cycle, has all curves satisfying $d_i\geq t_i+t'_i-2$ for all $i$, except that one allows that either
\begin {enumerate}
\item one curve satisfies $d=t+t'-3$, or 
\item two curves, separated by a (possibly empty) string of rational curves, satisfy $d=t+t'-3$.
\end{enumerate}
Then $h^1(S(-E'))=0$, $h^1(S)=h^1(S\otimes \mathcal O_E),$ and the Rational Conjecture is satisfied.
\begin{proof}  Combining Lemmas \ref{e2} and \ref{e3} and Remark \ref{r1}, we conclude that $H^1(S(-E')\otimes \mathcal O_E))=0.$  It suffices to show that $H^1(S(-E')\otimes \mathcal O_{nE})=0$ for all $n\geq 2$.   For each $n\geq 1$, one proceeds inductively from the divisor $nE$ to $(n+1)E$ via $nE+F$, for a judiciously chosen $F\geq 0$ which is effective and reduced.  Consider for a curve $E_i$ not contained in $F$ the sequence $$0\rightarrow S(-E'-nE-F)\otimes \mathcal O_{E_i}\rightarrow S(-E')\otimes \mathcal O_{nE+F+E_i}\rightarrow S(-E')\otimes \mathcal O_{nE+F}\rightarrow 0.$$  The requirement for the induction is that $H^1$ of the first term is $0$.  One has that $S(-E'-nE)\otimes \mathcal O_{E_i}$ equals
\begin{enumerate}
\item $\mathcal O_{E_i}(d_i-t'_i+n(d_i-t_i))\oplus \mathcal O_{E_i}(2-t_i+d_i-t'_i+n(d_i-t_i))$, if $t_i>1$
\item  $\mathcal O_{E_i}(-1+n(d_i-1))\oplus \mathcal O_{E_i}(n(d_i-1))$, if $t_i=1$.
\end{enumerate}
So, $H^1$ of the twist with $\mathcal O_{E_i}(-F)$ equals $0$ as long as 
\begin{enumerate}
\item $F\cdot E_i \leq 3-t_i-t'_i+d_i+n(d_i-t_i)$, if $t_i \geq 2$
\item $F\cdot E_i \leq n(d_i-1)$, if $t_i=1$.
\end{enumerate}

One way to proceed is to first choose some $E_0$ and go from $nE$ to $nE+E_0$, and then step by step add a curve adjacent to what has already been chosen (i.e., go from $nE+F$ to $nE+F+E_i$ if $F\cdot E_i=1).$  Given that $d_i\geq t_i$ and $d_i\geq t_i+t'_i-3$, then this procedure will work starting with any $E_0$ unless there is a curve with $$3-t_i-t'_i+d_i+n(d_i-t_i)=0,$$ i.e. $d_i=t_i$ and $t'_i=3$.  If there is only one such curve, then let $E_0$ be that curve; then the above procedure will get one from $nE$ to $nE+E_0$ and then on to $(n+1)E$, and the desired vanishing of $H^1(S(-E'))$ holds.

Now suppose there are two curves $E_0$ and $E_0'$ with $d=t+t'-3$.  If the path between them contains a curve $E_1$ with
$d_1>t_1+t'_1-2$, then start with $F=E_0+E'_0$, and successively adds the curves between $E_0$ and $E'_0$, up to $E_1$.  At this point, the new $F$ will satisfy $F\cdot E_1=2$, but now the inequality of $(1)$ is satisfied, and the induction may proceed as before.  

We are left with the case that the only curves (if any) in between $E_0$ and $E_0'$ satisfy $d=t+t'-2$, i.e. $d=t=t'=2$.  In that case, let $F$ be the sum of $E_0$, $E'_0$ and all the curves in between.  To proceed in the induction from $nE$ to $nE+F$, it suffices to show that $H^1(S(-E'-nE)\otimes \mathcal O_F))=0$.  But by assumption, every curve $E_i$ in $F$ satisfies $E\cdot E_i=0$, hence $\mathcal O_F(-nE)\cong \mathcal O_F.$  But $H^1(S(-E')\otimes \mathcal O_F)$ is a quotient of $H^1(S(-E')\otimes \mathcal O_E)$, which is $0$ as already mentioned.
\end{proof}
\end{theorem}
\begin{example}  Example \ref{ex1} is taut, i.e., $H^1(S)=0$.
\end{example}
\begin{example}\label{7.11} Assuming $b\geq 3$, the graph below is rational if either $a\geq 3$, or $a=2$ and $b\geq 5$:
$$
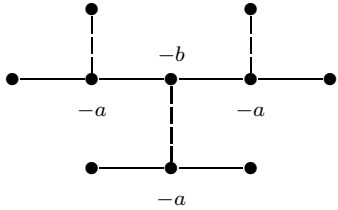
\xymatrix@R=6pt@C=24pt@M=0pt@W=0pt@H=0pt{
\\
&&&\overtag{\bullet}{}{6pt}&
&\overtag{\bullet}{}{6pt}&&&\\
&&&\lineto[u]&&\lineto[u]&\\
&&&\lineto[u]&&\lineto[u]&\\
&&\undertag{\bullet}{}{2pt}\lineto[r]
&\undertag{\bullet}{-a}{1pt}\lineto[r]\lineto[u]
&\overtag{\bullet}{-b}{10pt}\lineto[r]&\undertag{\bullet}{-a}{1pt}\lineto[r]\lineto[u]
&\undertag{\bullet}{}{2pt}
&\\
&&&&\lineto[u]&\\
&&&&\lineto[u]&&\\
&&&&\lineto[u]&\\
&&
&\overtag{\bullet}{}{8pt}\lineto[r]
&\undertag{\bullet}{-a}{1pt}\lineto[r]\lineto[u]
&\undertag{\bullet}{}{6pt}&&&\\
&&&
&&&
&&&\\
}
$$ 
\medskip

We ask whether $H^1(S(-E'))=0$, i.e. $h^1(S)=2.$  This is true in the following cases:
\begin{enumerate}
\item For $a,b\geq 4$, by the ``Easy Vanishing Theorem'' \ref{rt}.
\item For $a\geq 3$ and $b\geq 4$, by Corollary \ref{c1}.
\item For $a\geq 3$ and $b\geq 3$ by Theorem \ref{t3}.
\item For $a=2$ and $b\geq 7$ by Proposition \ref{van} applied to $D$ equals $E$ plus the $3$ outer nodes, concluding $H^1(S(-E'))\cong H^1(S(-E')\otimes \mathcal O_{D-E'})$, and noting the last term is a quotient of $H^1(S(-E')\otimes \mathcal O_E)$, which is $0$ by Theorem \ref{t1}.
\end{enumerate}
Our methods can not handle the case $a=2$ and $b=5$ or $6$, though $d_i\geq t_i+t'_i-2$, all $i$.
\end{example}

 \section{Results in characteristic $p$}
   Throughout this section we consider rational surface singularities in characteristic $p>0$.   We analyze earlier proofs to find sufficient conditions for the same calculations of $h^1(S)$ to hold.  Arguments that use Riemann-Roch (e.g., easy vanishing theorems, the Euler characteristic of $S\otimes \mathcal O_{E'}$) remain valid.  So we restate Theorem \ref{rt}: 
   \begin{theorem}\label{rtp}  If $d_i\geq 2t_i-2$ for all $i$, then in all characteristics $h^1(S)=h^1(S\otimes  \mathcal O_E)$.
   \end{theorem}   
    One needs  to revisit the calculation of $h^0(S\otimes  \mathcal O_{E'})$.  
   \begin{lemma}\label{l1p} Exclude cyclic quotients, cusps, and simple ellliptic singularities.
   \begin{enumerate}
   \item If the graph $\Gamma$ is star-shaped, then $h^0(S\otimes \mathcal O_{E'})=1.$
   \item Suppose $\Gamma$ contains two stars connected by a chain of rational curves, the determinant $n$ of whose intersection matrix is not divisible by $p$.   Then $h^0(S\otimes \mathcal O_{E'})=0.$  
   \end{enumerate}
   \begin{proof}  The previous proofs of Lemmas \ref{l4} and \ref{l5} are valid in characteristic $p$, given that the hypothesis in $(2)$ implies that the equations for the $B_i$'s admit only the $0$ solution.
   \end{proof}
   \end{lemma}

 Examining the proof of Theorem \ref{t1}, the arguments involving Type II cones involved some inequalities; we avoid that case by restricting to the case of reduced fundamental cycle (cf. Remark \ref{r1}).  The Riemann-Roch argument of Lemma \ref{s2} is valid in characteristic $p$, so a string of length at least two starting from an end-curve is still of Type I.  
  
  On the other hand, according to Lemma \ref{s3} an end-curve's contribution $(1,-d)$ becomes $(1,0)$ when $p$ divides $d$.  In this case, its neighbor receives $(0,1)$ in the corresponding slot.  For each curve $E_i$, define $\bar{t_i}$ to be the number of adjacent end-curves whose degrees are divisible by $p$.  Note $\bar{t_i}\leq t_i-t'_i.$ The situation is clarified by the following analogue of Lemma \ref{beg}.
  \begin{lemma}\label{lp}Suppose $\mathcal C(E_0,P_0)$ is a Round $2$ cone (so $t'=1$), and assume $d\geq t$.  Then 
 \begin{enumerate}
\item $d<t+\bar{t}-2$ implies  $\Phi_{E_0,P_0}$ is not surjective, hence $\Phi_E$ is not surjective
\item $d=t+\bar{t}-2$ implies $\mathcal C$ is Type II, except that the contribution at $P_0$ is of the form $(0,1)$.
\item $d>t+\bar{t}-2$ implies $\mathcal C$ is Type I.
\end{enumerate}
\begin{proof}  Since $d\geq t$, as in the proof of Lemma \ref{beg} we can always produce contributions of $(0,1)$ at each of the $t$ points of $E_0$.  But for $\bar{t}$ of the points, we need a contribution of the type $(1,*)$.  This requires having an $A(x)$ which vanishes at all but one of these points; this means that $$d-t+1\geq \bar{t}-1,$$ or $d\geq t+\bar{t}-2.$  The proof should now be clear.
\end{proof}
 \end{lemma} 
 One avoids the new type II condition in the last result via the inequality $$d\geq t+t'+\bar{t} -2.$$  These are exactly the conditions one needs to generalize all previous results from characteristic $0$, since the induction involves only end-curves and Type I curves.
 \begin{theorem}\label{t4}  Suppose a rational singularity, with reduced fundamental cycle, has all curves satisfying $d_i\geq t_i+t'_i+\bar{t_i}-2$ for all $i$, except that one allows that either
\begin {enumerate}
\item one curve satisfies $d=t+t'+\bar{t}-3$, or 
\item two curves, separated by a (possibly empty) string of rational curves, satisfy $d=t+t'+\bar{t}-3$.
\end{enumerate}
Then $h^1(S(-E'))=0$, so $h^1(S)=h^1(S\otimes \mathcal O_E)$.
\begin{proof} To prove first that $H^1(S(-E')\otimes \mathcal O_{E})=0$, consider the inductive step of Lemma \ref{midround} for a curve $E_0$.  There are no Type I curves, so the only slots that need filling are from the $t-t'$ end-curves plus the curve at $P_0$.  Contributions of the type $(0,1)$ are handled by choosing $A(x)=0$ and various $B(x)$ to vanish at all but one of these points.  This can happen because $d-t'\geq t-t'$.  We also need contributions of type $(1,0)$ at $P_0$ and the $\bar{t}$ points.  This requires choosing various $A(x)$ to vanish at all but one of these points; but $d-t-t'+2\geq \bar{t}$ by hypothesis, so this can happen.   A similar argument handles the terminal situation, except that choosing $A(x)$ to vanish at all but one of the $\bar{t}$ points requires only the weaker condition that $d-t-t'+2\geq \bar{t}-1$ (the situation in $(1)$ above.)  
     Moving to the situation of $(2)$, assume that the induction has led to a cone $\mathcal C$ with the weaker condition $d=t+t'+\bar{t}-3$.  Then the argument above shows the cone has Type II, except that the contribution at $P_0$ is $(0,1)$.  Now consider the case that one has two such cones, separated by a string of rational curves.  Then the argument in Lemma \ref{e3} works exactly as before.  Therefore, in all cases of the Theorem one has $H^1(S(-E')\otimes \mathcal O_{E})=0$. To conclude that $H^1(S(-E'))=0$, the proof of Theorem \ref{t3} is valid in all characteristics.
\end{proof}
 \end{theorem}
  \begin{corollary}  Suppose the graph is star-shaped, and the central curve satisfies $d\geq t$.  Then $h^1(S(-E')\otimes \mathcal O_E)=0$ if and only if $$d\geq t+t'+\bar{t}-3.$$
 \begin{proof}  At the node, there are $t'$ Type I curves.  Letting $A(x)=0$, one can choose various $B(x)$ of degree $d-t'$ to vanish at all but one of the $t-t'$ end-curve intersection points. If $\bar{t}=0$, then automatically $\Phi_E$ is surjective, so for vanishing of $H^1$ one only needs that $d\geq t+t'-3$ (obvious converse to Lemma \ref{Cok}).
 If $\bar{t}>0$, then to separate out those points requires finding various $A(x)$ of degree at least $\bar{t}-1$; this means $d-t-t'+2\geq \bar{t}-1.$ 
 \end{proof}
 \end{corollary}
 \begin{remark} The inequality in the Corollary is automatic if $d\geq 2t-3$, so the case $d=t=3$ is covered.  However, if $d=t=\bar{t}=4$ (a degree $4$ central curve plus $4$ end-curves whose degrees are divisible by $p$),  then the cohomology group does not vanish, and there are ``extra" equisingular deformation beyond those arising from the cross-ratio of the $4$ intersection points on the central curve.
 \end{remark}
  \begin{example}\label{pex}  Consider a star-shaped graph whose central curve has $d=4$, and each of whose $4$ branches consists of a single $-p$ curve:

 $$
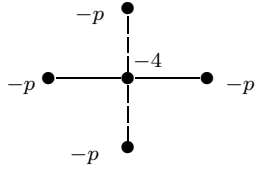
\xymatrix@R=6pt@C=24pt@M=0pt@W=0pt@H=0pt{
\\
&&&\lefttag{\bullet}{-p}{6pt}&
&&&&\\
&&&\lineto[u]&&&\\
&&&\lineto[u]\righttag{}{-4}{2pt}&&&\\
&&\lefttag{\bullet}{-p}{2pt}\lineto[r]
&\undertag{\bullet}{}{2pt}\lineto[r]\lineto[u]
&\righttag{\bullet}{-p}{10pt}&
&\\
&&&\lineto[u]&\\
&&&\lineto[u]&&\\
&&
&\lefttag{\bullet}{-p}{8pt}\lineto[u]
&&&\\
&&&
&&&
&&&\\
}
$$ 
 The argument in the proof of Lemma \ref{lp} shows that in characteristic $p$, one has $h^1(S(-E')\otimes \mathcal O_E)\neq 0$, hence $h^1(S)\geq 2.$  
  \end{example} 
    \begin{remark}  Recall that in characteristic $0$, the  ``hard'' vanishing theorem $H^1_E(S)=0$ implies that the space of equisingular deformations of a rational resolution inject into a smooth subspace of the base space of the semi-universal deformation of the singularity \cite{vt}.  That vanishing result need no longer be true in characteristic $p$, although it is not known whether the equisingular deformations still inject into the base-space; the non-vanishing may simply reflect the failure to lift vector fields from the singularity to the MGR.  We note, however, that if the fundamental cycle is reduced and at least one non-end curve has $d_i>t_i$, then the vanishing theorem still holds (by \cite{vt}, (2.16)).
  
  \end{remark}

 \section{Taut singularities in characteristic $p$ with reduced fundamental cycle}
      We shall use the following criterion to determine whether a graph $\Gamma$ is taut.
 \begin{theorem}(\cite{lauf},(3.9); \cite{sch})  A graph $\Gamma$ is taut if and only if for every singularity with graph $\Gamma$, on the MGR one has $H^1(S)=0.$
 \begin{proof}  Note that if $H^1(S)\neq 0$, then there would be a non-trivial smooth equisingular family of resolutions; but tautness implies one has a unique singularity, hence a unique resolution.
 
     Next suppose $H^1(S)=0$ for every resolution.  For an effective cycle $Z$ on a resolution, there is an easily verified surjection $S\rightarrow \Theta_Z$.  Thus, $H^1(\Theta_Z)=0$.   Laufer takes a graph $\Gamma$ which is ``potentially taut'' (i.e., all $t_i\leq 3$) and a formal sum $Z=\sum n_iE_i$, converting it into a ``plumbing scheme'' $P=P_Z$; this is an actual exceptional divisor on a resolution of a specific singularity with graph $\Gamma$.   The requisite characteristic $p$ construction is similar, done in Section $3$ of \cite{sch}.  The authors show that if $Z$ is sufficiently big, then $H^1(P,\Theta_P)=0$ implies tautness (\cite{sch}, Proposition 3.16).  The point is that a combinatorially equivalent divisor on another resolution can be connected to $P$ by a connected family (actually a more general result is proved by Laufer in \cite{lauf}, Theorem 3.2).
     \end{proof}
 \end{theorem}  
 
          If $h^1(S)=0$, then so is $h^1(S\otimes \mathcal O_E)$.  By Proposition \ref{p6}, except for the excluded cases, taut singularities are star-shaped with $3$ ends or are not star-shaped and have $4$ ends.  One thus considers also $\Gamma$ a chain of rational curves (a cyclic quotient) or a cycle of rational curves (a ``cusp" singularity).   In the following, every vertex is allowed to have any degree $\geq 2$ unless otherwise specified.
      
      \begin{theorem}\label{tautp} The following are the taut singularities in characteristic $p$ with reduced fundamental cycle:
      \begin{enumerate}
      \item For all $p$,$$
\xymatrix@R=12pt@C=24pt@M=0pt@W=0pt@H=0pt{
&\undertag{\bullet}{}{3pt}\dashto[r]&\dashto[r]&
\undertag{\bullet}{}{3pt}&\\
}
$$
\item For $d\geq 3$ and all $p$,  
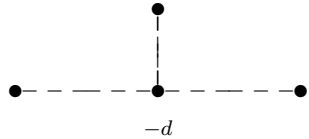
$$
\xymatrix@R=12pt@C=24pt@M=0pt@W=0pt@H=0pt{
&&{\bullet}\dashto[dd]\\
&&\dashto[u]\\
{\bullet}\dashto[r]&&
\undertag{\bullet}{-d}{3pt}\dashto[l]\dashto[u]\dashto[r]&\dashto[r]&
\undertag{\bullet}{}{3pt}\\
&&&&\\
}
$$
\item For $a,b\geq 3$ and $p$ not dividing the determinant of the string of curves between the nodes, $$
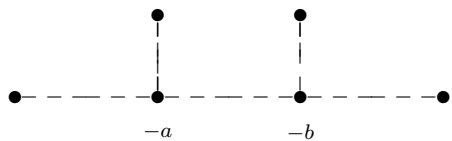
\xymatrix@R=12pt@C=24pt@M=0pt@W=0pt@H=0pt{
&&{\bullet}\dashto[dd]&&{\bullet}\dashto[dd]\\
&&\dashto[u]&&\dashto[u]\\
{\bullet}\dashto[r]&&
\undertag{\bullet}{-a}{2pt}\dashto[l]\dashto[u]\dashto[r]&\dashto[r]&
\undertag{\bullet}{-b}{2pt}\dashto[r]&\dashto[r]&{\bullet}\\
&&&&\\
}
$$ 
\item For $p$ not dividing the determinant of the intersection matrix of the cusp,  $$
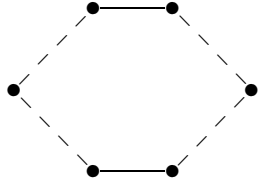
\xymatrix@R=12pt@C=24pt@M=0pt@W=0pt@H=0pt{
&&{\bullet}\dashto[ddl]\lineto[r]&{\bullet}\dashto[ddr]\\
&&&&\\
&{\bullet}\dashto[ddr]&
&&{\bullet}\dashto[ddl]\\
&&&&\\
&&{\bullet}\lineto[r]&{\bullet}\\
&&&&\\
}
$$ 
 \end{enumerate}
\begin{proof}  Except for restrictions on the prime $p$, the above list includes all taut singularities in characteristic $0$ with reduced fundamental cycle.

Theorem \ref{rtp} says that $h^1(S)=h^1(S\otimes \mathcal O_E)$ in all cases above except $(2)$ for $d=3$ and $(3)$ for $a$ or $b$ equal to $3$.  But those cases are covered by Theorem \ref{t4}.

    In cases $(2)$ and $(3)$, Proposition \ref{p6} gives the value of $h^1(S\otimes \mathcal O_E)$ as long as 
    $h^0(S\otimes \mathcal O_E)$ is as it was in characteristic $0$.  By Lemma \ref{l1p}, we conclude that for $(2)$ and $(3)$, one has $h^1(S)=0$; these are indeed taut. 
    
    For $(3)$, we show that if $p$ does divide the determinant, then $h^1(S\otimes \mathcal O_E)\neq 0$.  For, there then exist non-trivial solutions $A_i$ and $B_i$ (in the notation of the proof of Lemma \ref{l4}), hence $h^0(S\otimes  \mathcal O_F)=1$.  Now go from $F$ to $E'$ arguing as in the proof to show $h^0(S\otimes \mathcal O_{F})=h^0(S\otimes  \mathcal O_{E'})=1$.
   The first line of the proof of Proposition \ref{p6} gives that $\chi(S\otimes \mathcal O_{E'})=0$, hence $h^1(S\otimes \mathcal O_{E'})=1$, and so $h^1(S)=1$.
   
     For $(1)$,  start with any curve in the string and proceed as in Lemmas \ref{l1} and \ref{l4} to conclude that $h^0(S\otimes \mathcal O_E)=4$; the Euler characteristic is also $4$, so $h^1(S)=0$.
     
     Finally, for the cusps of $(4)$, the Euler characteristic of $S\otimes \mathcal O_E$ is $0$, so $h^1(S\otimes \mathcal O_E)=h^0(S\otimes \mathcal O_E)$.  A calculation similar to that in the proof of Lemma \ref{l5} gives equations for the coefficients $A_i,B_i$ of sections of $H^0(S\otimes \mathcal O_{E_i})$; they reduce to homogeneous equations in the $B_i$ whose determinant is that of the intersection matrix of the graph.  Thus, there is a non-trivial solution if and only if $p$ divides this determinant.

\end{proof} 
 \end{theorem}
 \begin{remark}  The graphs pictured for the cusps of $(4)$  show more than one curve, but the same argument applies when the minimal resolution is a nodal curve whose degree is not divisible by $p$.
 \end{remark}
 \begin{remark} Lee and Nakayama have already proved \cite{LN} that the cyclic quotient singularities in $(1)$ are taut in all characteristics.  Other cases of tautness in characteristic $p$ have been proved by Y. Tanaka \cite{tan} and F. Sch\"{u}ller \cite{sch}.
 \end{remark}
 \begin{remark} Of course, there are many star-shaped rational graphs (e.g., $E_8$) which are taut in characteristic $0$, but not in certain positive characteristic (e.g., \cite{ma}).  Necessarily, those must have $d$ equal $2$.  \end{remark}

\bigskip

\end{document}